\input amstex

\def\b1{\text{\bf 1}}

\def\CA{{\Cal A}}
\def\CB{{\Cal B}}
\def\CC{{\Cal C}}
\def\CD{{\Cal D}}
\def\CF{{\Cal F}}

\def\CJ{{\Cal J}}
\def\CG{{\Cal G}}
\def\CL{{\Cal L}}
\def\CM{{\Cal M}}

\def\CO{{\Cal O}}
\def\CP{{\Cal P}}

\def\CS{{\Cal S}}

\def\CV{{\Cal V}}

\def\gr{\text{gr}}

\def\Hom{\text{Hom}}

\def\et{\text{\'et}}
\def\h{\text{h}}
\def\hotimes{\widehat{\otimes}}

\def\#{\,\check{}}

\def\fm{{\frak m}}

\def\fD{{\frak D}}

\def\ft{{\frak t}}

\def\fk{{\frak k}}
\def\fa{{\frak{a}}}

\def\fs{{\frak s}}

\def\nc{{\text{nc}}}

\def\ss{{\text{ss}}}

\def\id{\text{id}}

\def\Ker{\text{Ker}}

\def\Spec{\text{Spec}}

\def\Z{{\text{Zar}}}
\def\dR{\text{dR}}

\def\BdR{\text{B}_{\text{dR}}}
\def\AdR{\text{A}_{\text{dR}}}

\def\limleft{\mathop{\vtop{\ialign{##\crcr
  \hfil\rm lim\hfil\crcr
  \noalign{\nointerlineskip}\leftarrowfill\crcr
  \noalign{\nointerlineskip}\crcr}}}}

\def\hra{\hookrightarrow}
\def\iso{\buildrel\sim\over\rightarrow} 

\def\lra{\longrightarrow}


\documentstyle{jams}
\issueinfo{25}{3}{Xxxx}{2012}
\NoBlackBoxes


\topmatter  \title p-adic periods and derived de Rham cohomology \endtitle  \author A.~Beilinson \endauthor\address Department of Mathematics, University of Chicago, Chicago, IL 60637 \endaddress \email sasha\@math.uchicago.edu \endemail  \thanks The author was  supported in part by NSF grant DMS-1001660. \endthanks    \subjclassyear{2010}   \subjclass Primary 14F30, 14F40; Secondary 14F20  \endsubjclass    \keywords p-adic periods, derived de Rham cohomology, h-topology, alterations \endkeywords      \date February 22, 2011 and, in revised form, November 16, 2011, and January 5, 2012   \enddate    
 \dedicatory To Irene  \enddedicatory  \endtopmatter
 
\document
 
\head Introduction. \endhead

For a smooth variety $X$ over a base field of characteristic 0 we have its algebraic de Rham cohomology $H^\cdot_{\dR}(X):=H^\cdot (X_\Z ,\Omega^\cdot_X )$; for nonsmooth $X$, one defines 
 $H^\cdot_{\dR}(X)$ using cohomological descent  as in Deligne \cite{D}. If the base field is $\Bbb C$,
 then one has the Betti cohomology $H^\cdot_{\text{B}}(X):=H^\cdot (X_{\text{cl}}, \Bbb Q )$ and  a canonical period isomorphism (``integration of algebraic differential forms over topological cycles")
$$\rho  : H^\cdot_{\dR}(X)\iso H^\cdot_{\text{B}} (X )\otimes\Bbb C  \tag 0.1$$    compatible with the Gal$(\Bbb C/\Bbb R )$-conjugation. To define $\rho$, 
consider the analytic de Rham cohomology $H^{\cdot}_{\dR }(X_{\text{an}})$. 
There are evident maps $$H^\cdot_{\dR}(X)\buildrel{\alpha}\over\to H^\cdot_{\dR} (X_{\text{an}})\buildrel{\beta}\over\leftarrow H^\cdot_{\text{B}}(X)\otimes\Bbb C . \tag 0.2$$  Then $\beta$ is an isomorphism due to the Poincar\'e lemma, and $\rho :=\beta^{-1}\alpha$ (the fact that $\rho$ is an isomorphism was established by Grothendieck \cite{Gr}). 

Suppose our base field is an algebraic closure $\bar{K}$ of a p-adic field $K$ (say, $K=\Bbb Q_p$). The role of $H^\cdot_{\text{B}}(X)$ is now played  by the p-adic \'etale cohomology $H^\cdot_{\et} (X,\Bbb Q_p )$, and Fontaine  conjectured in \cite{F1}\footnote{The assumption of loc.~cit.~that $X$ is  proper and smooth is redundant.} the existence of a natural p-adic period isomorphism $$\rho : H^\cdot_{\dR}(X)\otimes_{\bar{K}}\BdR \iso H^\cdot_{\et} (X,\Bbb Q_p )\otimes\BdR . \tag 0.3$$ Here $\BdR$ is Fontaine's p-adic periods field (\cite{F1}, \cite{F3}). Recall that  it  is a complete discretely-valued field whose ring of integers $\BdR^+$ contains $\bar{K}$, the residue field $\BdR^+/\fm_{\dR} $ is Tate's field $\Bbb C_p$, the cotangent line $\fm_{\dR} /\fm_{\dR}^2$ is the Tate twist  $\Bbb C_p (1)$. Both sides of 
 (0.3) carry natural filtrations (coming from the filtration of $\BdR$ by powers of $\fm_{\dR}$ and the Hodge-Deligne filtration on $H^\cdot_{\dR}(X)$), and $\rho$ is compatible with them and with the Gal$(\bar{K} /K )$-conjugation. Moreover, as was envisioned by Fontaine and Jannsen  (\cite{F4}, \cite{J}), the matrix coefficients of $\rho$  lie in the  subring $\bar{K} \text{B}_{\text{st}}$  of $\BdR$,\footnote{This assertion is true for a general reason, since, by Berger's theorem \cite{Ber},  de Rham Galois modules are potentially log crystalline.}  and $\rho$ is compatible with the extra symmetries of  log crystalline cohomology.

The p-adic period map  was defined in three  different ways in  works of, respectively, Faltings, Niziol, and Tsuji (with  prior crucial input of Bloch, Fontaine, Hyodo, Kato, Kurihara, and Messing; the nonproper setting was treated by Yamashita \cite{Y}), see  \cite{Fa1}, \cite{Fa2},  \cite{N1}, \cite{N2}, \cite{Ts1}, \cite{Ts2}; the  three  $\rho$'s coincide by \cite{N3}.

In the article we give another construction of $\rho$ which  is fairly direct and has the same flavor as  the classical picture (0.2). The  tools are derived de Rham cohomology of Illusie \cite{Ill2} Ch.~VIII and de Jong's alterations. The companion paper \cite{B}
treats the Fontaine-Jannsen side of the story;  another approach was developed by Bhatt \cite{Bh2}. It would be very interesting to see if these methods can help to understand the Riemann-Hilbert correspondence in the p-adic setting.

\enspace

An outline of the construction: First we realize 
B$^+_{\text{dR}}$   as the ring of de Rham p-adic constants in the sense of {\it derived} algebraic geometry. Namely, let 
$\AdR$ be the derived de Rham algebra  $L\Omega^\cdot\!\hat{_{O}}\!{}_{{}_{\bar{K}}  /O_K}$ completed with respect to the Hodge filtration $F^\cdot$, see \cite{Ill2} Ch.~VIII (2.1.3.3);   here  $O_K$, $O_{\bar{K}}$ are the rings of integers in $K$, $\bar{K}$.  
Now   $\BdR^+ $ identifies canonically with $\AdR \hotimes \Bbb Q_p   $, where $\hotimes$ is  the  derived completed tensor product, so that $\fm_{\text{dR}}^i \iso F^i  \AdR \hotimes \Bbb Q_p $. This fact  was observed independently by  Fargues \cite{Far}.

Let $\CV ar_F$ be the category of varieties over a field $F$, and
 $\CV ar^{\nc}_F$ be the category of regular $F$-varieties $U$ equipped with a regular compactification $\bar{U}$ with normal crossings divisor at infinity. 
 As follows from de Jong's theorem \cite{dJ1}, the  forgetful functor $\CV ar^{\nc}_F \to \CV ar_F$, $(U,\bar{U})\mapsto U$, makes  $\CV ar^{\nc}_{F}$  a base for the h-topology on $\CV ar_{F}$, so h-sheaves on $\CV ar_{F}$ are the same as sheaves on $\CV ar^{\nc}_{F}$ for the induced topology. 
 For $F=\bar{K}$ as above, there is a finer category 
$\CV ar^{\ss}_{\bar{K}}$ of {\it ss-pairs} $(V,\bar{V})$, i.e., smooth $\bar{K}$-varieties $V$  equipped with a semi-stable compactification $\bar{V}$ (that includes compactification in the arithmetic direction). 
Again by de Jong  \cite{dJ1},   $\CV ar^{\ss}_{\bar{K}}$
 is a base for the h-topology on  $\CV ar_{\bar{K}}$.

Consider the presheaf on $\CV ar^{\ss}_{\bar{K}}$ which assigns to $(V,\bar{V})$   the derived de Rham algebra with log singularities $R\Gamma (\bar{V},L\Omega^\cdot\!\hat{_{(V}}\!{}_{,\bar{V})/O_K})$  (see \cite{Ol}). Its h-sheafification $\CA^\natural_{\dR}$ is an h-sheaf 
of filtered dg algebras on $\CV ar_{\bar{K} }$ that  
 contains the constant subsheaf $\AdR$. The key {\it p-adic Poincar\'e lemma} says  that {\it the map   $\AdR\otimes^L \Bbb Z /p^n \to \CA^\natural_{\dR}\otimes^L \Bbb Z /p^n$ is a filtered quasi-isomorphism}. It comes from the next assertion: The h-sheafification of the presheaf
$(V,\bar{V})\mapsto  H^b (\bar{V},\Omega^{a}_{(V,\bar{V})/O_{\bar{K}}})$,  where $\Omega^{a}_{(V,\bar{V})/O_{\bar{K}}}$ is the usual locally free $\CO_{\bar{V}}$-module of forms with log singularities, is an h-sheaf of $\Bbb Q$-vector spaces
 for $(a,b)\neq (0,0)$.  The case $a=0$  is essentially theorem 8.0.1 from Bhatt's thesis \cite{Bh1}; the general result is obtained by a similar method (which uses coverings of families of stable curves that come from the multiplication by $p$ isogeny of the generalized Jacobians). 

 Set $R\Gamma^\natural_{\!\dR}(X ):= R\Gamma (X_{ \h},\CA^\natural_{\dR} )$; this is the {\it arithmetic de Rham complex} of $X$. By the above, $ H^\cdot (R\Gamma^\natural_{\!\dR}(X )\hotimes\Bbb Q_p )$ is a $\BdR^+$-algebra. One has a diagram 
$$  H^\cdot_{\dR}(X) \buildrel\alpha\over\to H^\cdot (R\Gamma^\natural_{\!\dR}(X )\hotimes\Bbb Q_p )
 \buildrel\beta\over\leftarrow H^\cdot_{\et}  (X,  \Bbb Q_p )\otimes \BdR^+ , \tag 0.4$$ where $\alpha$ is the composition 
$ H^\cdot_{\dR}(X)\iso H^\cdot R\Gamma^\natural_{\!\dR}(X )\otimes \Bbb Q \to H^\cdot (R\Gamma^\natural_{\!\dR}(X ) \hotimes\Bbb Q_p )$ and $\beta$ is the $\BdR^+$-linear extension of the evident map (which comes from the embeddings $\Bbb Z/p^n \to \CA^\natural_{\dR}\otimes^L \Bbb Z/p^n$ and the fact that the h-topology is stronger than the \'etale one).
Since the \'etale and h-cohomology with  torsion coefficients coincide, the Poincar\'e lemma implies that $\beta$ is an isomorphism. Now the p-adic period map $\rho$ is the $\BdR$-linear extension of $\beta^{-1}\alpha$.
An explicit computation for $X=\Bbb G_m$ followed by usual tricks of the trade shows that $\rho$ is a filtered isomorphism.

 \enspace
 
The work was prompted by a discussion with Kobi Kremnizer (and  watching the moon over a meadow in Oxford). I am very grateful to him, and to Bhargav Bhatt, Volodya Drinfeld, 
John Francis,
Luc Illusie, David Kazhdan,  Jacob Lurie,  Madhav Nori, and Misha Temkin for explanations, questions, and correction of  mistakes. 

\head 1.  A derived de Rham construction of B$_{\text{dR}}$. \endhead

\subhead{\rm 1.1}\endsubhead {\it The derived p-adic completion.} Throughout the article we use (not too heavily) E$_\infty$ algebras, for which we refer to, say, 
\cite{HS}.\footnote{There E$_\infty$ algebras are called ``May algebras".} Recall that
 E$_\infty$ algebras are dg algebras whose  product  is commutative and associative up to coherent higher homotopies (more formally, E$_\infty$ algebras are dg algebras for a resolution of the commutative algebra operad). A key fact: for any commutative (more generally, E$_\infty$)  cosimplicial dg algebra the corresponding total complex is naturally an E$_\infty$ algebra. Thus the homotopy limit of a diagram of  E$_\infty$ algebras is an E$_\infty$ algebra.

For a projective system of complexes of abelian groups $\ldots  \buildrel{\phi_2}\over\lra C_2 \buildrel{\phi_1}\over\lra C_1$, one has holim $C_n =\CC one (\id -\phi :\Pi C_n \to \Pi C_n )[-1]$, where $\phi ((c_n ))= (\phi_n (c_{n+1}) )$.  There is an embedding $ \limleft C_n =\Ker (\id- \phi )\hra$  holim $ C_n $.
 If all $\phi_n$'s are surjective, then $\id -\phi$ is surjective, hence $\hra$ is a quasi-isomorphism. So holim, being an exact functor, is the right derived functor of $\limleft$. 

If $C_\cdot$ is a projective system of dg algebras, then holim $C_n$ is naturally a dg algebra (and the above embedding is an embedding of algebras); if the $C_n$ are commutative (or, more generally, E$_\infty$) algebras, then holim $ C_n$ is an E$_\infty$ algebra.

Let $p$ be a prime. Consider the projective system of commutative dg algebras
 $C_n :=\CC one (\Bbb Z \buildrel{p^n}\over\lra \Bbb Z )$. It is a resolution of the projective system $\ldots \to \Bbb Z/p^2 \to \Bbb Z/p $, so $\Bbb Z_p^\flat :=\text{holim}\, C_n$ is an E$_\infty$ algebra with $H^0  \Bbb Z_p^\flat =\Bbb Z_p$ and acyclic in nonzero degrees. Set $\Bbb Q_p^\flat :=  \Bbb Z_p^\flat \otimes \Bbb Q$.
 For any complex $F$ of abelian groups set $$F\hotimes \Bbb Z_p := \text{holim}\, (F\otimes C_n ), \quad F\hotimes \Bbb Q_p := (F\hotimes \Bbb Z_p )\otimes\Bbb Q . \tag 1.1.1$$
These are dg $\Bbb Z_p^\flat$ and $\Bbb Q_p^\flat$-modules, so their cohomology groups are $\Bbb Z_p$- and $\Bbb Q_p$-modules, and $F\mapsto   F\hotimes \Bbb Z_p ,F\hotimes \Bbb Q_p$ are exact functors.
If $F$ is an (E$_\infty$) dg algebra, then so are $F\hotimes \Bbb Z_p $ and $F\hotimes \Bbb Q_p$.

\remark{ Remark} One has an evident projective system $ F_{p^n}[1]\to F\otimes C_n \to F/p^n F $  of exact triangles; applying holim, we get a canonical exact triangle holim$ F_{p^n}[1] \to F\hotimes\Bbb Z_p \to \text{holim} (F/p^n F)$. Let  $\hat{F}:=\limleft F/p^n F$ be the p-adic completion of $F$ and $T_p F:= \limleft F_{p^n}$ be the Tate module of $F$. By above, we have a quasi-isomorphism $\hat{F}\iso \text{holim}( F/p^n F)$. Thus if $F$ has no p-torsion, then $
F\hotimes\Bbb Z_p \iso \hat{F}$. Similarly,  if all components of $F$ are p-divisible, then one has  quasi-isomorphisms $T_p F \iso \text{holim} F_{p^n}$ and $T_p F [1] \iso   F\hotimes\Bbb Z_p$. We see that $\cdot\hotimes\Bbb Z_p$ is the left derived functor of the p-adic completion functor and the right derived functor of $T_p [1]$.\endremark

\remark{Example} For a scheme $X$, 
 its \'etale $\Bbb Z_p$- and $\Bbb Q_p$-cohomology are $R\Gamma_{\!\et}(X,\Bbb Z_p ):= \text{holim} \, R\Gamma (    X_{\et},\Bbb Z /p^n )=R\Gamma (X_{\et},\Bbb Z)\hotimes\Bbb Z_p$,  $R\Gamma_{\!\et}(X,\Bbb Q_p ):= R\Gamma (X_{\et},\Bbb Z)\hotimes\Bbb Q_p$.\footnote{Since $ R\Gamma (X_{\et},\Bbb Z/p^n ) =  R\Gamma (X_{\et},\Bbb Z)\otimes^L \Bbb Z/p^n$.} \endremark

\subhead{\rm 1.2}\endsubhead {\it The derived de Rham algebra.} For a morphism of commutative rings $A\to B$ we denote by $\Omega^\cdot_{B/A}$  the relative de Rham complex of $B$ over $A$. This is a commutative dg $A$-algebra with  components $\Omega^i_{B/A}=\Lambda^i_B  \Omega_{B/A}$, where $\Omega_{B/A}$ is the $B$-module of relative K\"ahler differentials; it carries a  ring filtration $F^n = \Omega^{\ge n}_{B/A}$. 

We will use the $F$-completed version 
$L\Omega^\cdot\!\!\hat{_{B/A}}$ of Illusie's
 derived de Rham algebra  defined in  \cite{Ill2} Ch.~VIII, (2.1.3.3). 
To construct it, consider the canonical simplicial resolution $P_\cdot = P_A (B)_\cdot$ of $B$ from \cite{Ill1}  Ch.~I,  (1.5.5.6). This is a simplicial commutative $A$-algebra such that each $P_i$ is a polynomial $A$-algebra. The de Rham complexes $\Omega^\cdot_{P_\cdot /A}$ form a simplicial filtered commutative dg $A$-algebra, so the corresponding total complex $L\Omega^\cdot_{B/A} $
 is  a  filtered commutative dg $A$-algebra (see \cite{Ill1} Ch.~I, 3.1.3). Now $L\Omega^\cdot\!\!\hat{_{B/A}}$
 is its completion with respect to the filtration $F^\cdot$.  {\it Here ``completion" is  understood as mere projective system of quotients modulo $F^i$}. One has a natural identification $\gr^i_F L\Omega^\cdot\!\!\hat{_{B/A}}\iso (L\Lambda^i_B (\text{L}_{B/A}))[-i]$ compatible with the product; here 
 L$_{B/A}:=\Omega_{P_\cdot /A}\otimes_{P_\cdot}B$ is the relative cotangent complex  and $L\Lambda^i_B$ is the nonabelian left derived functor of the exterior power functor   (see Ch.~II and I of  \cite{Ill1}).   For $A$-flat $B$'s, the construction is compatible with base change. It is  compatible with direct limits.  If in the above definition we replace $P_\cdot$ by any
simplicial $A$-algebra resolution  of $B$ whose terms are polynomial $A$-algebras, then the output is naturally quasi-isomorphic to $L\Omega^\cdot\!\!\hat{_{B/A}}$.

\enspace

The next lemma is a particular case of \cite{Ill1} Ch.~I,   4.3.2.1(ii). For a flat $B$-module $T$ we denote by $
B\langle T\rangle^\cdot
$ its divided powers symmetric algebra.

\proclaim{Lemma}  The complex $  L \Lambda_B^i (T[1])$ is acyclic off degree $-i$. There is a canonical isomorphism of graded $B$-algebras compatible with base change $$ H^{-\cdot} L \Lambda_B^\cdot (T[1]) \iso B\langle T\rangle^\cdot
. \tag 1.2.1$$\endproclaim

\subhead{\rm1.3}\endsubhead Let  $K$ be a p-adic field, i.e., a complete discretely-valued field of characteristic zero with perfect residue field  $k$ of characteristic $p>0$, $\bar{K}$ be an algebraic closure of $K$, and $O_K$, $O_{\bar{K}}$ be rings of integers in $K$, $\bar{K}$. 
Let $K_0 \subset K$ be the field of fractions of the Witt vectors $W(k)=O_{K_0}$,  and let $\fa$ be the fractional ideal in $\bar{K}$ generated by $p^{-\frac{1}{p-1}}\fD^{-1}_{K/K_0}$, where $\fD_{K/K_0}$ is the different. For an $O_K$-algebra $B$ we often write $\Omega_B := \Omega_{B/O_K}$, $L\Omega^\cdot\!\hat{_B}:= L\Omega^\cdot\!\!\hat{_{B/}}{}_{O_K}$, L$_B =\text{L}_{B/O_K}$, etc.

The next key result is due to Fontaine   \cite{F2} Th 1; we include a proof for completeness sake. Consider 
the map $\mu_{p^\infty} \subset O_{\bar{K}}^\times \, \buildrel{d\log}\over\lra \, \Omega_{O_{\bar{K}}}$ and its  $O_{\bar{K}}$-linear extension
$$(\bar{K}/O_{\bar{K}})(1)= O_{\bar{K}}\otimes \mu_{p^\infty} \to \Omega_{O_{\bar{K}}}. \tag 1.3.1$$

\proclaim{ Theorem}  One has L$_{O_{\bar{K}}}\iso\Omega_{O_{\bar{K}}} $, and 
  (1.3.1) 
 is surjective with kernel $(\fa /O_{\bar{K}})(1)$.\endproclaim
 
 \demo{Proof}  If $K'/K$ is a finite extension, then $O_{K'}/O_K$ is a complete intersection. So, if
$\pi$ is a generator of $O_{K'}/O_K$, $f(t)$ its minimal polynomial,
then L$_{O_{K'}}$ is the cone of multiplication by $f'(\pi )$ endomorphism of $O_{K'}$;
hence L$_{O_{K'}} \iso \Omega_{O_{K'}}$. Passing to the limit, we get the first assertion. Let us prove the second one. 

 (i) By the above, $\Omega_{O_{K'}}\simeq O_{K'}/\fD_{K'/K}$. If $K''/K'$ is another finite extension, then the standard exact triangle of the cotangent complexes reduces to a short exact sequence $0\to O_{K''}\otimes_{O_{K'}}\Omega_{O_{K'}/O_K}\to\Omega_{O_{K''}/O_K}\to
\Omega_{O_{K''}/O_{K'}}\to 0$. 

(ii) Replacing $K'$, $K$ by $K$, $K_0$ and passing to the limit, we get a short exact sequence $0\to O_{\bar{K}}\otimes_{O_{K}}\Omega_{O_{K}/O_{K_0}}\to\Omega_{O_{\bar{K}}/O_{K_0}}\to\Omega_{O_{\bar{K}}/O_{K}}\to 0$. Thus it suffices to prove the theorem for $K=K_0$, which we now assume.

(iii)  Set $T:=\Ker ( (\bar{K}/O_{\bar{K}})(1) \to \Omega )$, $F:=K (\mu_p )$.
The set of $O_{\bar{K}}$-submodules of $ (\bar{K}/O_{\bar{K}})(1)$ is  totally ordered by inclusion. Thus, since $O_{\bar{K}}\otimes_{O_{F}}  \Omega_{O_{F}} \subset \Omega_{O_{\bar{K}}}$ is a nonzero $O_{\bar{K}}$-module generated by $d\log (\mu_p )$, one has $T\subset (p^{-1}O_{\bar{K}}/O_{\bar{K}})(1)=O_{\bar{K}}\otimes\mu_p $. Since 
$ \Omega_{O_{F}}$ is isomorphic to $O_{F}/p^{1-\frac{1}{p-1}}O_{F}$, one has $T= (p^{-\frac{1}{p-1}}O_{\bar{K}}/O_{\bar{K}})(1)$.

(iv) It remains to prove surjectivity of $(\bar{K}/O_{\bar{K}})(1) \to \Omega_{O_{\bar{K}}}$. Let $K'\subset\bar{K}$ be any finite extension of $K$;  
we want to check that $\Omega_{O_{K'}}\subset \Omega_{O_{\bar{K}}}$ lies in $O_{\bar{K}}d\log (\mu_{p^\infty})$. Suppose $p^n $ kills $\Omega_{O_{K'}}$. Let us show that $\Omega_{O_{K'}}\subset O_{\bar{K}}d\log (\mu_{p^{n+1}})$. Set $K'' :=K'(\mu_{p^{n+1}})$. The set of $O_{K''}$-submodules of $\Omega_{O_{K''}}$ is  totally ordered. Thus, since $p^n d\log (\mu_{p^{n+1}})\neq 0$ by (iii), $\Omega_{O_{K'}}$ lies in $ 
 O_{K''}d\log (\mu_{p^{n+1}})$, q.e.d.      \qed
  \enddemo

\subhead{\rm1.4}\endsubhead  {\it For a complex $P$ acyclic in  degrees $\neq 0$,  we often write $P$ instead of $H^0 P$.}

\enspace

Consider the filtered commutative  dg $O_K$-algebra $\AdR =\text{A}_{\dR\,\bar{K}/K}  := L\Omega^\cdot\!\!\hat{_{O_{\bar{K}}}}{}_{/O_K}$ and the corresponding filtered E$_\infty$  $O_K$-algebra
 $\AdR \hotimes \Bbb Z_p$ (see 1.1). Let us describe the graded $O_{\bar{K}}$-algebras $\gr^\cdot_F \AdR$, $\gr^\cdot_F \AdR\hotimes \Bbb Z_p$.

\proclaim{Proposition}  (i) The complexes $\gr_F^i \AdR \hotimes \Bbb Z_p $ are acyclic in nonzero degrees,  and there is a canonical isomorphism of graded algebras $$  \gr^\cdot_F \AdR \hotimes \Bbb Z_p \iso   \hat{O}_{\bar{K}}\langle\hat{\fa}(1)\rangle^\cdot .
 \tag 1.4.1$$
\newline
(ii) One has $\gr^0_F \AdR =\AdR /F^1 = O_{\bar{K}}$, and the complexes $\gr^i_F \AdR$ for $i>0$ are acyclic in degrees $\neq 1$. There are natural isomorphisms of $O_{\bar{K}}$-modules
$$\Omega^{\langle i \rangle} :=H^1  \gr^i_F \AdR    \iso (\bar{K}/i!^{-1}\fa^i )(i)= (\Bbb Q_p /\Bbb Z_p ) \otimes i!^{-1} \hat{\fa}^i (i). \tag 1.4.2$$ \endproclaim

\demo{Proof} (i) By the theorem in 1.3, one has
L$_{O_{\bar{K}}/O_K}\iso\Omega_{O_{\bar{K}}}\iso
(\bar{K}/  \fa )(1)= ( \Bbb Q_p /\Bbb Z_p ) \otimes \fa (1)$; hence
 $\gr^i_F \AdR \iso L\Lambda^i_{\Bbb Z} (\Bbb Q_p / \Bbb Z_p )[-i] \otimes \fa^i (i)$. 
One has $L\Lambda^i_{\Bbb Z} ( \Bbb Q_p / \Bbb Z_p )\otimes^L (\Bbb Z/p^n )=L\Lambda^i_{\Bbb Z /p^n} (( \Bbb Q_p / \Bbb Z_p )\otimes^L (\Bbb Z/p^n ))=  L\Lambda^i_{\Bbb Z/p^n} ((\Bbb Z /p^n )[1])$, which identifies with  
$i!^{-1} (\Bbb Z/p^n )[i]$ in a way compatible with the product by (1.2.1). Therefore 
$\gr^\cdot_F \AdR  \otimes^L (\Bbb Z/p^n )   \iso \Bbb Z/p^n \langle (\fa /p^n \fa )(1)\rangle^\cdot$, which yields (1.4.1). \newline (ii) follows from (i) by the
next observation  (applied to $C=\gr^i_F \AdR$, with (1.4.2) defined by the condition that $T_p$(1.4.2)=(1.4.1)): If a complex $C$ of abelian groups has p-torsion cohomology and 
$H^{\neq 0} (C\otimes^L \Bbb Z/p )=0$, then $H^1 C$ is p-divisible and $H^{\neq 1} C=0$.\footnote{Use the fact that every complex of abelian groups splits, i.e., is quasi-isomorphic to a complex with zero differential.} 
\quad\qed  \enddemo

\subhead{\rm1.5}\endsubhead By 1.4(i), the algebras $(\AdR /F^i )\hotimes \Bbb Z_p $, hence $(\AdR /F^i )\hotimes \Bbb Q_p $, are acyclic in non-zero degree.   By loc.cit.,  
$(\AdR /F^{i+1} )\hotimes \Bbb Q_p $ is an i-truncated dvr with residue field $\Bbb C_p := \hat{O}_{\bar{K}}\otimes \Bbb Q$, so $ \AdR \hotimes \Bbb Q_p :=\limleft (\AdR /F^i )\hotimes \Bbb Q_p$ is a dvr. Let $\fm_{\dR}$ be its maximal ideal;  (1.4.1) yields a canonical identification $\fm_{\dR}/\fm_{\dR}^2 =\gr_F^1  \AdR \hotimes \Bbb Q_p \iso \Bbb C_p (1)$.

\proclaim {Proposition} There is a canonical ring isomorphism of filtered rings $$ u_{\Bbb Q} :\BdR^+ \iso  \AdR \hotimes \Bbb Q_p  . \tag 1.5.1$$ \endproclaim

\demo{Proof} The ring $(\AdR /F^{i+1} )\hotimes \Bbb Z_p $ is an infinitesimal p-adic $O_K$-thickening of $\hat{O}_K $ $= (\AdR /F^{1} )\hotimes \Bbb Z_p $ of order $\le i$ (see \cite{F3} 1.1). Let A$_{\text{inf}}/F^{i+1}$ be the universal thickening (\cite{F3} 1.3); we have a canonical map $u_i : \text{A}_{\text{inf}}/F^{i+1}\to (\AdR /F^{i+1} )\hotimes \Bbb Z_p$. Since 
 $\BdR^+ / F^{i+1}:= ( \text{A}_{\text{inf}}/F^{i+1} )\otimes \Bbb Q$ is an $i$-truncated dvr and
$u_1$ is an isomorphism by \cite{F3} 1.4.3, 
  $u_{i\Bbb Q}:  \BdR^+ / F^{i+1} \iso (\AdR /F^{i+1} )\hotimes \Bbb Q_p $. Set $u_{\Bbb Q}:=\limleft
 u_{i\Bbb Q}$.  \quad\qed  \enddemo 

\remark {Remarks} (i) The map $\AdR\to \AdR/F^1 = O_{\bar{K}}$ yields an isomorphism
 $\AdR \otimes \Bbb Q \iso \bar{K}$. Thus the  morphism $\AdR \otimes \Bbb Q \to 
\AdR \hotimes \Bbb Q_p$ equals the usual embedding $\bar{K}\hra \BdR^+$.
\newline (ii) For a finite extension $K'/K$, $K'\subset\bar{K}$,  the evident map $A_{\dR\,\bar{K}/K}\to A_{\dR\,\bar{K}/K'}$ yields an isomorphism $A_{\dR\,\bar{K}/K} \hotimes \Bbb Q_p\iso A_{\dR\,\bar{K}/K'} \hotimes \Bbb Q_p$ compatible with (1.5.1).
\endremark

\subhead{\rm1.6}\endsubhead The next result, which
will not be used in the rest of the article, is a reinterpretation of Colmez's theorem \cite{Col}. It would be nice to find a simpler direct proof.

\proclaim {  Proposition}  The complexes  $\AdR /F^i$ are acyclic in nonzero degrees;  the  maps $ H^0 (\AdR/F^{i+1} )\to H^0 ( \AdR/F^{i} )$ are injective. Set $O^{\langle i\rangle}:= H^0 ( \AdR/F^{i +1})$; thus $O_{\bar{K}}= O^{\langle 0\rangle}\supset O^{\langle 1\rangle}\supset \ldots$ and 
  $(\AdR/F^{i+1} )\hotimes \Bbb Z_p $ is equal to the p-adic completion $\hat{O}^{\langle i\rangle}$ of $O^{\langle i\rangle}$. 
\endproclaim

\demo{Proof}  By 1.4(ii),
 the  exact cohomology sequence for $0\to \gr^i_F \AdR \to \AdR /F^{i+1}$ $\to \AdR/F^{i}\to 0$ reduces to  $0 \to O^{\langle i\rangle}\to O^{\langle i-1\rangle}\buildrel{ d^{\langle i\rangle }  }\over\lra \Omega^{\langle i \rangle}
   \to H^1 (\AdR /F^{i+1})\to  H^1 (\AdR /F^{i})\to 0$. So $O_{\bar{K}}= O^{\langle 0\rangle}\supset O^{\langle 1\rangle}\supset \ldots$, and the  vanishing of $H^1 (\AdR /F^{i+1} )$ amounts to that of $H^1 (\AdR  /F^i )$ combined with surjectivity of  $d^{\langle i\rangle }:   O^{\langle i-1\rangle} \to  \Omega^{\langle i \rangle}$.  It remains to prove that all $d^{\langle i\rangle }$ are surjective. 
  
 Recall that Colmez   \cite{Col} considers  a sequence of subalgebras $O_{\bar{K}}=O^{(0)}\supset O^{(1)}\supset \ldots$ and derivations $d^{(i)}: O^{(i-1)}\to \Omega^{(i)}$ defined by induction: $d^{(i)}$ is a universal $O_K$-linear derivation with values in an $O_{\bar{K}}$-module, and $O^{(i)}:=\Ker\, d^{(i)}$. An induction by $i$ shows that  $O^{\langle i\rangle}\supset O^{(i)}$: Indeed,
 $\Omega^{\langle i \rangle}$ are $O_{\bar{K}}$-modules and $d^{\langle i\rangle }:   O^{\langle i-1\rangle} \to  \Omega^{\langle i \rangle}$ is a derivation; so, if  $O^{\langle i-1\rangle}\supset O^{(i-1)}$, then
 $d^{\langle i\rangle }|_{ O^{(i-1)}  }=a^{(i)}d^{(i)}$ for some  $O_{\bar{K}}$-linear map   $a^{(i)}: \Omega^{(i)}\to \Omega^{\langle i \rangle}$; thus $O^{\langle i\rangle}\supset O^{(i)}$.

Let $i$ be the smallest number such that $d^{\langle i\rangle }$ is not surjective. Since $E:=\Omega^{\langle i \rangle }/ d^{\langle i \rangle }(O^{\langle i -1\rangle})$ is p-torsion p-divisible, one has  $E\hotimes \Bbb Q_p = T_p E \otimes\Bbb Q\neq 0$. Applying $\cdot \hotimes\Bbb Z_p$ to the  exact triangle $ O^{\langle i \rangle }\to \AdR /F^{i+1}\to E[-1]$, we get a short exact sequence $0\to
\hat{O}^{\langle i \rangle }\to (\AdR /F^{i+1})\hotimes\Bbb Z_p \to T_p E \to 0$. By \cite{Col},  A$_{\text{inf}}/F^{i+1}=\hat{O}^{(i)}$. By universality, the map $u_{i}: \text{A}_{\text{inf}}/F^{i+1}\to   (\AdR /F^{i+1} )\hotimes \Bbb Z_p $ equals the composition $\hat{O}^{(i)}\to \hat{O}^{\langle i \rangle }\hra (\AdR /F^{i+1})\hotimes\Bbb Z_p$, so its composition with the projection onto $  T_p E$ vanishes. This cannot happen since $u_{i\Bbb Q}$ is an isomorphism (see 1.5), q.e.d.  
\quad\qed  
    \enddemo 

\head  2.   $h$-topology and  semi-stable compactifications. \endhead

\subhead{\rm 2.1}\endsubhead {\it A  topological digression.} The next proposition is a generalization  of \cite{V2} 4.1.

 Let $\CV$ be an essentially small site. As in \cite{V1}, we denote by  $\CV\,\,\widetilde{}\,$ the corresponding topos (the category of sheaves of sets on $\CV$).

 For us, a {\it base for $\CV$} is a pair $(\CB,\phi )$, where $\CB$ is an essentially small category and $\phi : \CB\to\CV$ is a {\it faithful} functor, that satisfies the next property: 
 \newline $(*)$ {\it For any $V\in\CV$ and a finite family of pairs $(B_\alpha ,f_\alpha )$,  $B_\alpha \in \CB$, 
  $f_\alpha : V\to \phi (B_\alpha )$,  there exists a set of objects  $B'_\beta \in\CB$ and
 a   covering family $\{ \phi (B'_\beta ) \to V \}$ such that every composition $\phi (B'_\beta ) \to V\to \phi (B_\alpha )$ lies in $\Hom (B'_\beta ,B_\alpha )\subset \Hom (\phi (B'_\beta ),\phi (B_\alpha ))$.}

\remark{Remarks} (i) Property
$(*)$ for  empty set of $(B_\alpha ,f_\alpha )$'s means that 
 every $V\in\CV$ has a covering by objects $\phi (B)$, $B\in\CB$. 
 If $\phi$ is fully faithful, then $(*)$ amounts to this assertion.\footnote{The proposition below in this situation amounts to
 \cite{V2} 4.1.}
\newline (ii) If $\CB$ admits finite products and $\phi$ commutes with finite products, then
it suffices to check $(*)$ for families $(B_\alpha ,f_\alpha )$ having $\le 1$ element.
 \newline (iii) In the general case, it suffices to check $(*)$ for families  $(B_\alpha ,f_\alpha )$ having $\le 2$ elements.
  \endremark

\enspace

Suppose $(\CB,\phi )$ is a base for $\CV$. Define a covering sieve in $\CB$ as a sieve whose $\phi$-image is a covering family in $\CV$. 

\proclaim {Proposition}  (i) Covering sieves in $\CB$ form a Grothendieck topology on $\CB$. \newline (ii) The functor $\phi : \CB\to \CV$ is continuous (see \cite{V2} 1.1). \newline (iii) $\phi$ yields an equivalence of  the toposes: one has $\CB\,\,\widetilde{}\,\, \iso \,\,\CV \,\,\widetilde{}$.
\endproclaim

We call the above topology on $\CB$ the {\it $\phi$-induced} topology.\footnote{The terminology is compatible with that of \cite{V2} 3.1. }

\demo{Proof} (i)
Let us check that  covering sieves in $\CB$ are stable with respect to pullback; the rest of the axioms from \cite{V1} 1.1 are evident. For a morphism  $g: B'\to B$  in $\CB$ and a covering sieve $\fs$ on $B$, let us find a covering family on $B'$ that belongs to the $g$-pullback of $\fs$.
The $\phi(g)$-pullback of $\phi (\fs )$ is a covering sieve in $\CV$, so
 there is a  covering family $\{\pi_\gamma : V_\gamma \to \phi (B')\}$ such that every composition
$V_\gamma \to \phi (B')\to\phi (B)$ can be factored as $V_\gamma \buildrel{g_\gamma}\over\lra \phi (B_{\gamma } ) \buildrel{\phi (p_\gamma )}\over\lra \phi (B)$, where $p_\gamma :  B_\gamma \to B$ belong to $\fs$. Applying $(*)$ to $V_\gamma$ and  $(B', \pi_\gamma )$, $(B_{\gamma },g_\gamma )$,
 we find a covering family $\{ \phi (B'_{\beta_{\gamma}}) \to V_\gamma \}$ as in $(*)$. 
The composite covering $\{ \phi (B'_{\beta_{\gamma}}) \to \phi (B')\}$   comes then from a covering family $\{ B'_{\beta_\gamma}\to B'\}$ in $\CB$
 which lies in the $g$-pullback of $\fs$.

(ii) We know that $\phi$ sends covering families to covering families, so it suffices to show that for any given $p_\alpha : B_\alpha \to B$ in $\CB$ and $f_\alpha : V\to \phi (B_\alpha )$,
$\alpha=1,2$, such that $\phi (p_1 ) f_1 =\phi (p_2 )f_2$ 
there is a covering  $\{\pi_\beta : V_\beta \to V\}$ and morphisms $\xi_{\alpha\beta}: B'_\beta \to B_\alpha$, $g_\beta : V_\beta \to \phi (B'_\beta )$ such that $p_1\xi_{1\beta}=p_2 \xi_{2\beta}$ and $\phi ( \xi_{\alpha\beta})g_\beta = f_\alpha \pi_\beta$. Such a datum (with $g_\beta$ the identity map) comes from $(*)$ applied to $V$ and $(B_1 , f_1 ), (B_2 ,f_2 )$.

(iii) By (ii), one has the usual adjoint functors between the categories of sheaves $(\phi^s ,\phi_s) : \CB\,\,\widetilde{} \leftrightarrows \CV\,\,\widetilde{}$. To prove that they are mutually inverse equivalences, we will check that for $\CF \in \CB\,\,\widetilde{}$ and $\CG \in\CV \,\,\widetilde{}\,$ the adjunction maps $a_{\CF}: \CF \to \phi_s \phi^s \CF$, $b_{\CG}: \phi^s\phi_s \CG \to\CG$ are isomorphisms.

Recall that $ \phi^s \CF =(\phi^\cdot \CF )\,\tilde{}$, where $\phi^\cdot$ is the pullback of presheaves and $\,\tilde{}\,$ is the sheafification functor. For $V\in\CV$ one has $(\phi^\cdot \CF )(V)=\text{colim}_{\CC (V)}\CF$, where $\CC (V)$ is the category
of pairs $(B,f)$,  $B\in\CB$, $f: V\to \phi (B)$, with $\Hom_{\CC (V)} ((B,f), (B',f')): = \{ g\in\Hom (B',B): \phi(g)f'=f\}$, and we set $\CF (B,f):=\CF (B)$.

(a) To show that $a_{\CF}$ is an isomorphism, we check that it is  injective and surjective:

{\it $a_{\CF}$ is injective}: Suppose we have $B\in \CB$ and $\xi_1 ,\xi_2 \in\CF (B)$ such that $a_{\CF}(\xi_1)=a_{\CF}(\xi_2 )$;
let us show that $\xi_i$ coincide. One has $a_{\CF}(\xi_i )\in (\phi_s \phi^s \CF )(B)= (\phi^s \CF )(\phi (B))$, and the equality means that there is a covering  $\{\pi_\gamma :  V_\gamma \to \phi (B)\}$ such that the images of $\xi_i$ in $(\phi^\cdot \CF )(V_\gamma )=    \text{colim}_{\CC (V_\gamma )}\CF $ coincide. Thus for some {\it finite} subdiagram $\CC (V_\gamma )' \subset \CC (V_\gamma )$ that contains $(B,\pi_\gamma )$ the images of $\xi_i$ in $\text{colim}_{\CC (V_\gamma )'}\CF$ coincide. Applying $(*)$ to $V_\gamma$ and pairs from $\CC (V_\gamma )'$, we get a covering  $\{ \phi (B'_{\beta_\gamma})\to  V_\gamma \}$ such that the image of $\CC (V_\gamma )'$ in each $\CC (\phi (B'_{\beta_\gamma}))$ comes from a diagram in $\CB^\circ /B'_{\beta_\gamma}$. The composite covering $\{ \phi (B'_{\beta_\gamma})\to \phi (B) \}$ comes then from a covering $\{ B'_{\beta_\gamma}\to B\}$ in $\CB$, and   the images of $\xi_i$ in $\CF (B'_{\beta_\gamma})$ coincide. Then $\xi_1 =\xi_2$
since $\CF$ is a sheaf,  q.e.d.

{\it $a_{\CF}$ is surjective}: For $B\in \CB$, $\chi \in     (\phi_s \phi^s \CF )(B)$ we look for a covering $\{ B'_\beta \to B \}$ in $\CB$ 
such that $\chi|_{B'_\beta }$ lies in the image of $\CF (B'_\beta )\to   (\phi_s \phi^s \CF )(B'_\beta )$. To find it, consider 
 $\chi$ as an element of $ (\phi^s \CF )(\phi(B))$. There is a covering $\{ \pi_\gamma : V_\gamma \to \phi (B) \}$ such that $\chi|_{V_\gamma}$ lies in the image of $(\phi^\cdot \CF )(V_\gamma )\to (\phi^s \CF )(V_\gamma )$, i.e., one has $f_\gamma : V_\gamma \to \phi (B_\gamma )$ 
 such that $\chi|_{V_\gamma}$ lies in the image of  the composition $\CF (B_\gamma )\to (\phi^s \CF )(\phi (B_\gamma ))\to  (\phi^s \CF )(V_\gamma )$, the second arrow comes from $f_\gamma$.
Applying $(*)$ to $V_\gamma$ and $(B,\pi_\gamma )$, $(B_\gamma, f_\gamma )$, we find a covering $\{ \phi (B'_{\beta_\gamma })\to V_\gamma \}$ as in $(*)$; the composite covering $\{ \phi (B'_{\beta_\gamma })\to \phi (B)\}$ comes then from a covering $\{ B'_{\beta_\gamma }\to B\}$
that satisfies the promised property.

(b) {\it $b_{\CG}$ is an isomorphism}: Since $\phi_s (b_{\CG})a_{\phi_s \CG }=\id_{\phi_s \CG}$ and we already know that $a_{\phi_s \CG}$ is an isomorphism, we see that $\phi_s (b_{\CG}): \phi_s \phi^s \phi_s (\CG )\to \phi_s \CG$ is an isomorphism. Thus $b_{\CG} (B) : 
 \phi^s \phi_s \CG (\phi (B))\to  \CG(\phi (B))$ is an isomorphism for every $B\in\CB$. Since every $V\in\CV$ admits a covering by objects $\phi (B)$, $B\in\CB$, this implies that $b_{\CG}$ is both injective and surjective, hence an isomorphism, q.e.d.  \quad\qed 
\enddemo

\remark{Exercises} (i) For any presheaf $\CJ$  on $\CV$ one has $\phi_s (\CJ\, \tilde{}\, )= (\phi_\cdot \CJ )\,\tilde{}$. \newline (ii) Suppose $(\CB,\phi )$ is a base for $\CV$ and $(\CB' ,\phi' )$ is a base for  the $\phi$-induced topology on $\CB$. Then $(\CB' ,\phi\phi' )$ is a base for  $\CV$. \endremark

\subhead{\rm2.2}\endsubhead For a field $K$, let  $\CV ar_K$ be the category of $K$-varieties, i.e., reduced separated $K$-schemes of finite type. We will consider categories $\CB$ formed by varieties equipped with appropriate compactifications, referred to as {\it pairs}:

(a) {\it Geometric setting:} Let $j: U\hra \bar{U}$ be an open embedding 
such that $\bar{U}$ is proper and $U$ is dense in $\bar{U}$. We call such a datum  {\it geometric pair over $K$}, or {\it geometric $K$-pair}, and denote it  by $(U,\bar{U})$. We say that $(U,\bar{U})$ is a regular normal crossings pair, {\it nc-pair} for short, if $\bar{U}$ is a regular scheme and $\bar{U}\smallsetminus U$ is a divisor with normal crossings in $\bar{U}$; it is a {\it strict} nc-pair if the irreducible components of $\bar{U}\smallsetminus U$ are regular.
A morphism $f: (U',\bar{U}')\to (U,\bar{U})$ of  pairs is a map $ \bar{U}'\to\bar{U}$ which sends $U'$ to $U$. We denote the category of geometric $K$-pairs by $\CV ar^{\text{c}}_K$; let  $\CV ar^{\nc}_K$ be the full subcategory of nc-pairs.

(b) {\it Arithmetic $K$-setting:} Suppose $K$ is a p-adic field as in 1.3. An {\it arithmetic pair over $K$}, a.k.a.~{\it arithmetic $K$-pair}, is an open embedding $j: U\hra \bar{U}$ with dense image of a $K$-variety $U$  into a reduced proper flat $O_K$-scheme $\bar{U}$. 

For such a $(U,\bar{U})$ we set  $O_{K_U} :=\Gamma (\bar{U}, \CO_{\bar{U}})$, $K_U := \Gamma (\bar{U}_K , \CO_{\bar{U}})$. Then $K_U$ is the product of several finite extensions of $K$ (labeled by the connected components of $\bar{U}_K$; if $\bar{U}$ is normal, then $O_{K_U}$ is the product of the corresponding rings of integers. {\it The closed fiber $\bar{U}_s$} of $\bar{U}$ is the union of fibers over the closed points of $O_{K_U}$.

We say that $(U,\bar{U})$  is a {\it  semi-stable pair}, or simply {\it ss-pair}, if   (i) $\bar{U}$  is a regular scheme, (ii) $\bar{U}\smallsetminus U$ is a divisor with  normal crossings on $\bar{U}$, and (iii) the closed fiber $\bar{U}_s$ is reduced. Our ss-pair is {\it strict}  if the irreducible components of $\bar{U}\smallsetminus U$ are regular. 
Arithmetic $K$-pairs form a category $\CV ar_K^{\text{cc}}$;
let $\CV ar^{\ss}_K $ be the full subcategory of ss-pairs. 

(c) {\it Arithmetic $\bar{K}$-setting:} For $K$  as in (b), let $\bar{K}$ be its algebraic closure. An {\it 
arithmetic pair over $\bar{K}$}, a.k.a.~{\it arithmetic $\bar{K}$-pair}, is an open embedding $j: V\hra \bar{V}$ with dense image of a $\bar{K}$-variety $V$  into a reduced proper flat $O_{\bar{K}}$-scheme $\bar{V}$. 
A connected $(V,\bar{V})$ is said to be {\it semi-stable}, a.k.a.~{\it ss-pair}, if  there exists an ss-pair $(U,\bar{U})$ over $K$  and a $\bar{K}$-point $\alpha :K_U \to \bar{K}$ (see (b)) such that  $(V, \bar{V})$ is isomorphic to $(U,\bar{U})_\alpha = (U_{\bar{K}},\bar{U}_{O_{\bar{K}}}):=
(U\otimes_{K_U}\!\bar{K},\, \bar{U}\otimes_{O_{K_U}}\!O_{\bar{K}})$. Then $\bar{V}$ is normal (say, by Serre's criterion). An arbitrary $(V,\bar{V})$ is semi-stable if such are all its connected components. Denote by $\CV ar^{\text{cc}}_{\bar{K}}$  the category of all arithmetic pairs over $\bar{K}$, and by $\CV ar^{\ss}_{\bar{K}}=\CV ar^{\ss}_{\!\bar{K}/K}$   its full subcategory of ss-pairs.

\remark{Remark} If $K'$ is a finite extension of $K$ contained in $\bar{K}$, then $\CV ar^{\ss}_{\!\bar{K}/K'}\subset \CV ar^{\ss}_{\!\bar{K}/K}$. For all the constructions below the difference between them is irrelevant.
\endremark

\enspace

These categories  are connected by  commutative diagrams of functors
$$\spreadmatrixlines{1\jot}
\matrix
\CV ar^{\text{cc}}_K &\to &
 \CV ar^{\text{c}}_{K} & \to &\CV ar_K , &&& \CV ar^{\text{cc}}_{\bar{K}}&\to&  \CV ar^{\text{c}}_{\bar{K}} &\to& \CV ar_{\bar{K}},   
 \\  
\uparrow  &&\uparrow&&&&& \uparrow&&\uparrow  \\
\CV ar^{\ss}_K &  \to & \CV ar^{\nc}_{K} &&&&& \CV ar^{\ss}_{\bar{K}}&\to&\CV ar^{\nc}_{\bar{K}}
    \endmatrix
\tag 2.2.1$$ where the vertical arrows are the fully faithful embeddings, and the upper horizontal lines are  faithful forgetful functors of passing to the generic fiber and  $(U,\bar{U})\mapsto U$. The $K$- and $\bar{K}$-settings are connected by base change functors
$$\spreadmatrixlines{1\jot}
\matrix
\CV ar^{\ss}_{\bar{K}}
 &\to &
 \CV ar^{\nc}_{\bar{K}} & \to &\CV ar_{\bar{K}}   \\  
\uparrow  &&\uparrow&& \uparrow  \\
\CV ar^{\ss}_K &  \to & \CV ar^{\nc}_{K} &\to& \CV ar_K .
    \endmatrix
\tag 2.2.2$$
Here the two right vertical arrows are the evident base change $\cdot \otimes_K \bar{K}$, and the left one assigns to a semi-stable $K$-pair $(U,\bar{U})$ the  disjoint sum of pairs $(U, \bar{U})_\alpha$ for all $\bar{K}$-points $\alpha : K_U \to \bar{K}$.

\subhead{\rm2.3}\endsubhead
A morphism $f: (V,\bar{V})\to (U,\bar{U})$ of pairs in either of the settings of 2.2 is called an {\it alteration} (of $(U,\bar{U})$) if $f^{-1}(U)=V$, the generic fibers of $f$  are zero-dimensional, and their union is dense in $V$. In setting (a),  $f$ is  a {\it (strict) nc-alteration} if $(V,\bar{V})$ is a (strict) nc-pair; in settings (b), (c), $f$ is a {\it (strict) ss-alteration} if $(V,\bar{V})$ is a (strict) ss-pair.

If $f$ is an alteration, then $f|_V :V\to U$ is proper and surjective; the composition of alterations is an alteration.

Here is a key result of de Jong  \cite{dJ1} 4.1, 6.5:
\proclaim{Theorem} Every geometric pair admits a strict nc-alteration. Every arithmetic pair, either over $K$ or over $\bar{K}$, admits a strict ss-alteration. The alterations can be chosen so that $\bar{V}$ is projective. 
 \quad\qed  \endproclaim

\remark{Remark} Our conventions slightly differ from de Jong's: he understands varieties to be irreducible and semi-stable $K$-pairs $(U,\bar{U})$ to have property $K_U = K$; his notation for $(U,\bar{U})$ is $(\bar{U}, Z)$,  $Z:=\bar{U}\smallsetminus U$.  \endremark 

\subhead{\rm2.4}\endsubhead For a field $K$,
 the {\it h-topology}  (see \cite{SV}) on   $\CV ar_K$  is generated by the pretopology whose coverings  are finite families of maps $\{ Y_\alpha \to X\}$ such that 
$Y:=\amalg Y_\alpha \to X$ is a universal topological epimorphism.\footnote{This means that a subset of $X$ is Zariski open if (and only if) its preimage in $Y$ is open, and the same is true after any base change.} It is stronger than the \'etale and proper topologies.\footnote{The latter is generated by a pretopology whose coverings are proper surjective maps.} We denote the $\h$-site by $\CV ar_{K\h}$; the h-site of $X$ is denoted by $X_{\h}$.

\remark{Exercise} Let $f: Y\to X$ be a morphism in $\CV ar_K$. \newline
(i) $f$ is an h-covering if (and only if) for every irreducible curve $C\subset X$ the base change $Y_{\tilde{C}}\to \tilde{C}$, where $\tilde{C}$ is the normalization of $C$, is an h-covering.\footnote{Hint: For an open  $V\subset Y$ its image in $X$ is constructible (EGA IV 1.8.4),  so to show that $f(V)$ is
open it suffices to check that for any curve $C\subset X$ the intersection $C\cap f(V)$ is open in $C$.} \newline
(ii) If $X$ is a regular curve, then $f$ is an h-covering if (and only if)  the closure of the generic fiber of $f$ maps onto $X$, or, equivalently, $f$  is a  covering for the flat topology. \endremark

\remark{Remark} By  \cite{SV} 10.4,  every  h-covering is a Zariski covering locally in proper topology.  Therefore  (see \cite{D}, \cite{SD}, or \cite{C})
 h-coverings are morphisms of universal cohomological descent for torsion \'etale sheaves; if $K=\Bbb C$, then h-coverings are morphisms of universal cohomological descent for arbitrary sheaves on the classical topology. In particular, for any h-hypercovering $Y_\cdot $ of $X$ and an abelian group $A$ the canonical map $R \Gamma(X_{\et} ,A)\to  R\Gamma(Y_{\cdot\, \et} ,A)$ (:= the total complex of the cosimplicial system of complexes $ R\Gamma(Y_{i\, \et} ,A)$)  is a quasi-isomorphism if $A$ is a torsion group. Passing to the limit, we see that $R\Gamma_{\!\et}(X,\Bbb Z_p )\iso R\Gamma_{\!\et}(Y_{\cdot } ,\Bbb Z_p )$.
If $K=\Bbb C$, then $R\Gamma (X_{\text{cl}},A)\iso R\Gamma (Y_{\cdot\text{cl} },A)$ for any $A$. Since  h-topology is stronger than the  \'etale one, we see  that 
  $R\Gamma (X_{\et} ,A)\iso R\Gamma (X_{\h},A)$ if $A$ is a torsion  group (see \cite{SV} 10.7 for a direct proof). \endremark

\subhead{\rm2.5}\endsubhead
Let $\phi$ be the forgetful functor $(U,\bar{U})\mapsto U$ on any of the categories $\CV ar^?$ in 2.2.

\proclaim{Proposition}  If  $\CV ar^? $ is either of the categories $ \CV ar^{\text{c}}_K $, $\CV ar^{\nc}_K$, $\CV ar^{\text{cc}}_K $, $\CV ar^{\ss}_K$, then $(\CV ar^? ,\phi )$ is a base  for $\CV ar_{K\h}$.  If $\CV ar^? $ is $ \CV ar^{\text{cc}}_{\bar{K}}$ or $ \CV ar^{\ss}_{\bar{K}}$, then $(\CV ar^? ,\phi )$ is a base for $\CV ar_{\bar{K}\h}$.
  \endproclaim

\demo{Proof} We consider the arithmetic $K$-setting, leaving the two other settings for the reader. Let us  show that  $(\CV ar^{\text{cc}}_K ,\phi)$ satisfies condition $(*)$ from 2.1. Our datum is a $K$-variety $V$ and a
finite collection of arithmetic $K$-pairs $(U_\alpha ,\bar{U}_\alpha )$ and maps $f_\alpha :V \to U_\alpha$. We need to find an h-covering $\pi : V'\to V$ and an arithmetic pair $(V',\bar{V}')$ such that $f_\alpha \pi$ extend to maps $(V',\bar{V}')\to (U_\alpha ,\bar{U}_\alpha )$. First we find an h-covering $V'\to V$ such that $V'$ sits in some arithmetic $K$-pair $(V',\tilde{V}')$: let $V'/V$ be a quasi-projective modification of $V$ provided by the Chow lemma, and take for $\tilde{V}'$ the closure of $V'$ in a projective space.\footnote{In fact, every $V$ sits in a $K$-pair due to Nagata's theorem.}  Then take  $\bar{V}'$ to be the closure of the image of $V'$ by the embedding $V'\hra \tilde{V}'\times \Pi \bar{U}_\alpha $, and we are done.

To show that $(\CV ar^{\ss}_K ,\phi )$ is a base for $\CV ar_{K\h}$, it suffices to check that
$(\CV ar^{\ss}_K ,\iota )$, where $\iota$ is the embedding $\CV ar^{\ss}_K \hra \CV ar^{\text{cc}}_K $, is a base for the $\phi$-induced topology on $\CV ar^{\text{cc}}_K$ (see Exercise (ii) in 2.1). Since $\iota$ is fully faithful, it suffices to check that for every  $(U,\bar{U})\in\CV ar_K^{\text{cc}}$ there exists a map $(U',\bar{U}')\to (U,\bar{U})$ such that $U'\to U$ is an h-covering and 
$(U',\bar{U}')$ is semi-stable. Such a datum is provided by the de Jong theorem in 2.3, and we are done. \qed
 \enddemo

We call the $\phi$-induced topology on either of the categories $\CV ar^?$ the {\it  h-topology}.

\remark{Remarks}  (i) Any h-covering of $(U,\bar{U})\in \CV ar^{\ss}_K$  has a refinement with terms of the same dimension as $U$ (indeed,  the same assertion  in $\CV ar_K$ is true by \cite{SV} 10.4; to pass to $\CV ar^{\ss}_K$, we apply the constructions from the proof above, and they preserve the dimension). 
\newline (ii) The proposition remains true if we replace the category of ss- or nc-pairs by its subcategory  of strict pairs $(U,\bar{U})$ with projective $\bar{U}$. 
\newline (iii) For any functor in (2.2.1) its source is a base for the h-topology of the target, and the induced topology on the source is the h-topology. 
\newline (iv) The  functors in (2.2.2) are continuous for the h-topologies.
\endremark

\subhead{\rm 2.6}\endsubhead 
 By 2.1 and 2.5, $\phi$ identifies  h-sheaves on $\CV ar_K$, resp.~ $\CV ar_{\bar{K}}$, with h-sheaves on $ \CV ar^{\text{c}}_K $, $\CV ar^{\nc}_K$, $\CV ar^{\text{cc}}_K $, $\CV ar^{\ss}_K$,
  resp.~$ \CV ar^{\text{cc}}_{\bar{K}}$, $ \CV ar^{\ss}_{\bar{K}}$. Thus we have the {\it h-localization} functors $$\CP\CS h (\CV ar_K^? )\to \CV ar\widetilde{_{K\h} } , \quad 
   \CP\CS h (\CV ar^{?}_{\bar{K}} )\to \CV ar\widetilde{_{\bar{K}\h} } \tag 2.6.1$$
 which assign to any presheaf $\CF$ on pairs the corresponding h-sheaf $\CF\,\,\tilde{}$ viewed as an h-sheaf on varieties. 
 
  \remark{Remark}  For any presheaf on $\CV ar^{\text{c}}_K$, $\CV ar^{\text{cc}}_K$ or $\CV ar^{\text{cc}}_{\bar{K}}$,  its h-sheafification coincides with h-sheafification of its restriction to resp.~$\CV ar^{\nc}_K$, $\CV ar^{\ss}_{K}$ or $\CV ar^{\ss}_{\bar{K}}$.
For a presheaf on $\CV ar^{\ss}_{\bar{K}/K}$, its h-sheafification is the same as h-sheafification of its restriction to $\CV ar^{\ss}_{\bar{K}/K'}$, where $K' \subset \bar{K}$ is any finite extension of $ K$ (see Remark in 2.2).
   \endremark

\head  3. The  p-adic period map. \endhead

\subhead{\rm 3.1}\endsubhead {\it The derived  de Rham algebra in logarithmic setting.} We refer to \cite{K1} for log scheme basics.
There are two (in general, nonequivalent) ways to define the cotangent complex for log schemes 
due, respectively, to  Gabber and Olsson, see \cite{Ol}.\footnote{In all situations we will consider the two versions coincide by \cite{Ol} 8.34.}
  Gabber's approach (\cite{Ol} \S 8) is more direct and precise;\footnote{It produces a true complex, while Olsson's construction yields a mere compatible datum of the canonical filtration truncations.}  we recall it briefly.

For a commutative ring $A$, a prelog structure on $A$ is a homomorphism of monoids $\alpha : L\to A$, where
$L$ is a commutative integral monoid (written multiplicatively) and $A$ is viewed as a monoid with respect to the product. 
Rings equipped with prelog structures  form a category in an evident way; denote its objects simply by $(A,L)$. For a fixed $(A,L)$,  let $\CC_{(A,L )}$ be the category of morphisms $(A,L )\to (B,M )$; we  denote such an object by $ (B,M )/  (A,L )$. Let $\Omega_{(B,M )/  (A,L )}$ be the $B$-module of relative K\"ahler log differentials: it is generated by $\Omega_{B/A}$ and elements $d\log m$, $m\in M$, subject to relations $d\log (m_1 m_2 )=d\log m_1 +d\log m_2 $, $\alpha (m)d\log m=d\alpha (m)$, and $d\log m=0$ if $m$ is in the image of $L$. The de Rham dg algebra of relative log forms $\Omega^\cdot_{(B,M )/  (A,L )}$ has components $\Omega^i_{(B,M )/  (A,L )}:=\Lambda^i_B  \Omega_{( B,M )/  (A,L )}$; elements $d\log m$ are degree 1 cycles. It carries the Hodge filtration $F^n =\Omega^{\ge n}_{(B,M)/  (A,L )}$.

A pair of sets $I$, $J$ yields a free object $P_{(A,L)} [I,J]$ in $ \CC_{(A,L )}$: the corresponding ring is a polynomial algebra $A[t_i , t_j]_{i\in I, j\in J}$, the monoid is $L \oplus  \Bbb N[I]$, where $\Bbb N[I]$ is the free monoid generated by $I$, and the structure map sends the generator $m_i$ of $\Bbb N[I]$, $i\in I$, to  $t_i $.  
The de Rham algebra $\Omega^\cdot_{P_{(A,L)} [I,J]/ (A,L)}$ is freely generated, as a graded commutative $A$-algebra, by elements $t_i$, $t_j$ of degree 0 and $d\log t_i :=d\log m_i$, $dt_j$ of degree 1, where $i\in I$, $j\in J$.

Every  $(B,M)/(A,L)\in\CC_{(A,L )}$ admits a canonical simplicial resolution $P_\cdot = P_{(A,L )}(B,M)_\cdot$. This is a simplicial object of $\CC_{(A,L )}$ 
 augmented over the object  $  (B,M)/(A,L)$ and such that every $P_i$ is a free object as above. Thus we have the simplicial dg algebra $\Omega^\cdot_{P_\cdot /(A,L)}$ filtered by the Hodge filtration $F$. Denote by $L\Omega^\cdot_{ (B,M)/(A,L)}$     the corresponding total complex,    $L\Omega^a_{ (B,M)/(A,L)}=\mathop\oplus\limits_{j-i=a} \Omega^j_{P_i  /(A,L)}$; this is a filtered commutative dg algebra. Let $L\Omega^\cdot\! \hat{_{ (B,}}{}_{M)/( A,L)}$ be its $F$-completion; as in 1.2, {\it we understand it as mere projective system of quotients $L\Omega^\cdot\! \hat{_{ (B,}}{}_{M)/( A,L)} /F^n$}. One has a natural quasi-isomorphism  of graded dg algebras  $\gr^*_F 
L\Omega^\cdot\! \hat{_{ (B,}}{}_{M)/(A,L)} \iso (L\Lambda^*_B (\text{L}_{(B,M)/(A,L)}))[-*].$ Here $\text{L}_{(B,M)/(A,L)} :=\Omega_{P_\cdot /(A,L)}\otimes_{P_\cdot}B$ is the relative log cotangent complex; it is acyclic in positive degrees, and $H^0 \text{L}_{( B,M)/(A,L)}=\Omega_{( B,M)/( A,L)}$.
The constructions are  compatible with direct limits.  If in the above definition we replace $P_\cdot$ by any free simplicial resolution  of $( B,M)/(A,L)$, then the output is naturally quasi-isomorphic to  $L\Omega^\cdot\! \hat{_{ (B,}}{}_{M)/(A,L)}$. The plain cotangent complex  and derived de Rham algebra for $B/A$ map naturally to logarithmic ones.

For any map $(X,\CM)\to (S,\CL )$ of integral log schemes, the above construction, being  \'etale sheafified, yields  the log cotangent complex L$_{(X,\CM)/(S,\CL)}$,
the derived log de Rham algebra  $L\Omega^\cdot_{(X,\CM)/(S,\CL )}$, and its $F$-completion $L\Omega^\cdot\! \hat{_{(X,}}{}_{\CM)/(S,\CL )}$, which are complexes of sheaves  on  $X_{\et}$. We use only the completed complex $L\Omega^{\cdot}\,\hat{}$.

\subhead{\rm 3.2}\endsubhead Let $(U,\bar{U})$ be a pair as in 2.2. We view
$\bar{U}$ as a log scheme with the usual integral log structure $\CO_{\bar{U}}\cap j_* \CO_U^\times  \to \CO_{\bar{U}}$; by abuse of notation, let us denote this log scheme again by $(U,\bar{U})$. 
Any morphism of pairs $(U,\bar{U})\to (V,\bar{V})$  is a morphism of log schemes, so we have the  relative log cotangent complex L$_{(U,\bar{U})/(V,\bar{V})}$, the derived log de Rham algebra $L\Omega^\cdot \!\hat{_{(U,}}{}_{\bar{U})/ (V,\bar{V})}$, etc., as above. There is a canonical morphism
$L\Omega^\cdot\! \hat{_{\bar{U}}}{}_{/\bar{V}}\to L \Omega^\cdot \!\hat{_{(U,}}{}_{\bar{U})/ (V,\bar{V})}$.
 We also have ``absolute" complexes: in the arithmetic $K$- or $\bar{K}$-setting, these are 
 L$_{(U,\bar{U})}:=\text{L}_{(U,\bar{U})/O_K}$, $L\Omega^\cdot \!\hat{_{(U,}}{}_{\bar{U})}:=  L\Omega^\cdot \!\hat{_{(U,}}{}_{\bar{U})/ O_K}$, where $O_K$ is considered with the trivial log structure $O_K^\times$; for the geometric $K$- or $\bar{K}$-setting,  replace $O_K$ by $K$, resp.~$\bar{K}$. 

\remark{Remark} For $(V,\bar{V})\in \CV ar^{\nc}_{\bar{K}}$ one has  $ L\Omega^\cdot \!\hat{_{(V,}}{}_{\bar{V})}\iso  \Omega^\cdot_{(V,\bar{V})}$. Hence for $(U,\bar{U})\in \CV ar^{\ss}_{\bar{K}}$ one has
$R\Gamma (\bar{U}, L\Omega^\cdot\!\hat{_{(U,}}{}_{\bar{U})})\otimes \Bbb Q \iso R\Gamma (\bar{U}_{\bar{K}}, \Omega^\cdot_{(U,\bar{U}_{\bar{K}})})$. Ditto for pairs over $K$.
\endremark

\enspace

Consider now the arithmetic $\bar{K}$-pair $\Spec (\bar{K},O_{\bar{K}}):=(\Spec\,\bar{K},\Spec\, O_{\bar{K}})$:

 \proclaim{Lemma} The cotangent complex
 $\text{L}_{(\bar{K},O_{\bar{K}})}$ is acyclic in nonzero degrees, and the canonical map $ \Omega_{O_{\bar{K}}}\to \Omega_{ (\bar{K},O_{\bar{K}})}:= H^0\text{L}_{(\bar{K},O_{\bar{K}})}$ is an isomorphism. Therefore the canonical map $\AdR:= L\Omega^\cdot\!\hat{_{O}}{}_{{}_{\bar{K}}  }\to L\Omega^\cdot\hat{_{(}}{}_{\bar{K},O_{\bar{K}})}$ is a filtered quasi-isomorphism.
 \endproclaim 

\demo{ Proof}  For a finite extension $K'$ of $K$ consider the log scheme $\Spec (K' ,O_{K'} ):=(\Spec\, K' ,\Spec\, O_{K'})$. It is a log complete intersection over $O_K$ (see \cite{Ol} 6.8). If $\pi$ is a generator of $O_{K'}/O_K$, $f(t)$ its minimal polynomial, then, by \cite{Ol} 6.9,      L$_{( K' ,O_{K'} )}$ is quasi-isomorphic to the cone of the multiplication by $f'(\pi)$ map $O_{K'}\to O_{K'}\subset \fm_{K'}^{-1}$. Thus
 L$_{( K' ,O_{K'} )}$ is acyclic in nonzero degrees, $\Omega_{( K' ,O_{K'} )}:= H^0 \text{L}_{( K' ,O_{K'} )}$ is a cyclic $O_{K'}$-module,
 and the canonical map $\Omega_{O_{K'}}\to \Omega_{( K' ,O_{K'} )}$ is an embedding with cokernel isomorphic to the residue field $O_{K'}/\fm_{K'}$. Now pass to the inductive limit, and use the fact that $\Omega_{O_{\bar{K}}}$ is p-divisible (see 1.3). \qed\enddemo

\subhead{\rm 3.3}\endsubhead Consider
 the presheaf $(U,\bar{U})\mapsto R\Gamma^\natural_{\!\dR}(U,\bar{U} ):=
 R\Gamma (\bar{U}, L\Omega^\cdot\!\hat{_{(U,}}{}_{\bar{U})})$ of filtered E$_{\infty}$ dg $O_K$-algebras on   $\CV ar^{\ss}_{\bar{K}}$. Denote by $\CA^\natural_{\dR}$ its  h-sheafification (2.6.1); this is an h-sheaf of filtered E$_\infty$  $O_K$-algebras on $\CV ar_{\bar{K}}$ (as above, we see it as the projective system of quotients modulo $F^i$). Since
  $\AdR  = \CA^\natural_{\dR} (\Spec \,\bar{K})$ by Lemma in 3.2,  $\AdR$, viewed as a constant filtered h-sheaf, maps into $\CA^\natural_{\dR}$. 

\proclaim { Theorem {\rm ($p$-adic Poincar\'e lemma)}} The maps $  \AdR\otimes^L \Bbb Z/p^n \to \CA^\natural_{\dR} \otimes^L \Bbb Z/p^n$ are filtered quasi-isomorphisms of h-sheaves on  $\CV ar_{\bar{K}}$. \endproclaim

For a proof, see \S 4. Assuming it, let us define the $p$-adic period map $\rho$.

\subhead{\rm 3.4}\endsubhead {\it The Hodge-Deligne filtration.}  For this subsection, $K$ is any field of characteristic 0.  Consider the presheaf 
$(V,\bar{V})\mapsto R\Gamma_{\!\dR}(V,\bar{V}):=
R\Gamma (\bar{V}, \Omega^\cdot_{(V,\bar{V})})$ of filtered E$_{\infty}$ dg $K$-algebras on   $\CV ar^{\nc}_{K}$. Let $\CA_{\dR}$ its  h-sheafification (2.6.1), which is an h-sheaf of filtered E$_\infty$  $K$-algebras on $\CV ar_{K}$ (viewed as the projective system of quotients modulo $F^i$). For any  $X\in \CV ar_{K}$ set   $$R\Gamma_{\!\dR}(X):= R\Gamma(X_{\h}, \CA_{\dR}). \tag 3.4.1$$ This is Deligne's de Rham complex of $X$ equipped with Deligne's Hodge filtration. 

\proclaim{Proposition} (i) For  $(V,\bar{V})\in \CV ar^{\nc}_{K}$ the canonical map $R\Gamma_{\!\dR}(V,\bar{V})\to R\Gamma_{\!\dR}(V)$ is a filtered quasi-isomorphism.
\newline (ii) The differential of $R\Gamma_{\!\dR}(X)$ is strictly compatible with the filtration.  $H^i_{\dR}(X):= H^i R\Gamma_{\!\dR}(X)$ are $K$-vector spaces of dimension equal to $
\dim H^i_{\et}(X_{\bar{K}},\Bbb Q_p )$. 
\newline (iii) For any smooth variety $X$  there is a canonical (nonfiltered) quasi-isomorphism $R\Gamma (X,\Omega^\cdot_X )\iso R\Gamma_{\!\dR}(X)$. \endproclaim

\demo{Proof} By  Lefschetz's principle, we can assume that $K=\Bbb C$. For  $(V,\bar{V})\in \CV ar^{\nc}_{\Bbb C}$ the maps $R\Gamma_{\!\dR}(V,\bar{V})
\to R\Gamma_{\!\dR}(V, \Omega^\cdot_V)\to R\Gamma (V_{\text{cl}},\Bbb C)$ are quasi-isomorphisms by \cite{Gr}. 
Thus for any h-hypercovering  $(Y_\cdot ,\bar{Y}_\cdot )/X$ of $X$ in $\CV ar^{\nc}_{K}$\footnote{Here we view $X$ as an h-sheaf on $\CV ar^{\nc}_{K}$, so  $(Y_\cdot ,\bar{Y}_\cdot )$ is a simplicial object of $\CV ar^{\nc}_{\bar{K}}$ equipped with an augmentation map $Y_\cdot \to X$ that makes $Y_{\cdot}$ an h-hypercovering of $X$.} the cohomological descent (see Remark in 2.4) yields a canonical quasi-isomorphism $
 R\Gamma (\bar{Y}_\cdot , \Omega^\cdot_{(Y_\cdot ,\bar{Y}_\cdot )})\iso 
 R\Gamma  (X_{\text{cl}},\Bbb C)$. If we equip $R\Gamma  (X_{\text{cl}},\Bbb C)$
 with the Hodge-Deligne filtration of mixed Hodge theory \cite{D}, then this
  is a filtered quasi-isomorphism. Therefore we have a canonical filtered quasi-isomorphism $R\Gamma_{\!\dR}(X)\iso R\Gamma (X_{\text{cl}},\Bbb C)$. Now (i) and the second assertion of (ii) are clear;
 the first assertion of (ii) follows from mixed Hodge theory. The quasi-isomorphism in (iii) is  $R\Gamma(X,\Omega^\cdot_X )\iso R\Gamma (Y_\cdot ,\Omega^\cdot_{Y_\cdot}) \buildrel\sim\over\leftarrow
R\Gamma (\bar{Y}_\cdot ,\Omega^\cdot_{(Y_.,\bar{Y}_\cdot )})  $, where the arrows are quasi-isomorphisms by the cohomological descent  (since $R\Gamma(X,\Omega^\cdot_X )\iso  R\Gamma (X_{\text{cl}},\Bbb C)$). \qed \enddemo

\subhead{\rm 3.5}\endsubhead We return to the setting of 3.3, so $K$ is our p-adic field. Let $X$ be any
variety  over $\bar{K}$. It yields a filtered  E$_\infty$  $O_K$-algebra
 $$ R\Gamma^\natural_{\!\dR}(X ):= R\Gamma (X_{\h},\CA^\natural_{\dR} ). \tag 3.5.1$$

Since $\AdR\otimes \Bbb Q =\bar{K}$ (see Remark (i) in 1.5),  $ R\Gamma^\natural_{\!\dR}(X )\otimes \Bbb Q$ is a $\bar{K}$-algebra. 
By Remark in 3.2, we have a filtered quasi-isomorphism of E$_\infty$  $\bar{K}$-algebras $$ R\Gamma^\natural_{\!\dR}(X )\otimes \Bbb Q \iso  R\Gamma_{\!\dR}(X). \tag 3.5.2$$

Let us compute $ R\Gamma^\natural_{\!\dR}(X )\hotimes \Bbb Z_p$. Consider the morphisms of filtered complexes  $R\Gamma (X_{\et} ,\Bbb Z)\otimes^L \AdR \iso R\Gamma (X_{\et} , \AdR)\to R\Gamma(X_{\h}, \AdR)\to R\Gamma(X_{\h}, \CA^\natural_{\dR} )= R\Gamma^\natural_{\!\dR}(X )$. After applying  $\cdot \otimes^L \Bbb Z/ p^n$, the arrows become filtered quasi-isomorphisms (the first one by Remark in 2.4, the second one by the Poincar\'e lemma in 3.3), so we get a filtered quasi-isomorphism $R\Gamma (X_{\et} ,\Bbb Z /p^n )\otimes^L \AdR \iso  R\Gamma^\natural_{\!\dR}(X )\otimes^L \Bbb Z/ p^n$. Since $R\Gamma_{\!\et}  (X,\Bbb Z_p )=\text{holim} R\Gamma (X_{\et} ,\Bbb Z/p^n )= R\Gamma (X_{\et} ,\Bbb Z)\hotimes \Bbb Z_p$ (see 1.1) is a perfect $\Bbb Z_p$-complex and $R\Gamma (X_{\et} ,\Bbb Z/p^n )= R\Gamma_{\!\et}  (X,\Bbb Z_p )\otimes^L_{\Bbb Z_p} \Bbb Z/p^n$,  one has, passing to the homotopy limit as in 1.1, $R\Gamma_{\!\et}  (X,\Bbb Z_p )\otimes^L_{\Bbb Z_p}(\AdR\hotimes \Bbb Z_p )\iso  R\Gamma^\natural_{\!\dR}(X )\hotimes\Bbb Z_p$. Tensoring by $\Bbb Q$, we get a filtered quasi-isomorphism of filtered  E$_\infty$  $\BdR^+$-algebras  (see (1.5.1))  $$ \beta :  R\Gamma_{\!\et}  (X,\Bbb Q_p )\otimes \BdR^+ \iso   R\Gamma^\natural_{\!\dR}(X ) \hotimes \Bbb Q_p .  \tag 3.5.3$$

 Let $\alpha : R\Gamma_{\!\dR}(X)\otimes_{\bar{K}} \BdR^+ \to
 R\Gamma^\natural_{\!\dR}(X ) \hotimes \Bbb Q_p$ be the $ \BdR^+$-linear extension of the
  composition $R\Gamma_{\! \dR}(X)\iso  R\Gamma^\natural_{\!\dR}(X )\otimes \Bbb Q \to  R\Gamma^\natural_{\!\dR}(X )\hotimes \Bbb Q_p$, where the first arrow is inverse to (3.5.2) and the second one comes from the canonical map $?\to ?\hotimes\Bbb Z_p$. We get a morphism of filtered
 E$_\infty$  $\BdR^+$-algebras $$\rho =\rho_{\dR}:=\beta^{-1} \alpha : R\Gamma_{\!\dR}(X)\otimes_{\bar{K}} \BdR^+ \to R\Gamma_{\!\et}  (X,\Bbb Q_p )\otimes_{\Bbb Q_p} \BdR^+ . \tag 3.5.4$$

 \remark{Remarks} (i) The Galois group Gal$(\bar{K}/K)$ acts on $\CV ar_{\bar{K}\text{h}}$ and on both sides of (3.5.4) by transport of structure, and $\rho^+$ is evidently compatible with this action. In particular, if
 $X$ is defined over $K$, i.e., $X=X_K \otimes_K \bar{K}$, then 
 $R\Gamma_{\!\dR}(X)= R\Gamma_{\!\dR}(X_K )\otimes_K {\bar{K}}$, and we can rewrite (3.5.4) as a Gal$(\bar{K}/K)$-equivariant morphism $$\rho  : R\Gamma_{\!\dR}(X_K )\otimes_{K} \BdR^+ \to R\Gamma_{\!\et}  (X,\Bbb Q_p )\otimes_{\Bbb Q_p} \BdR^+ . \tag 3.5.5$$ (ii) The map $\rho$ does not change if we replace $K$ by any its finite extension that contains in $\bar{K}$ (see Remark in 2.6).  \endremark

\proclaim{{\rm 3.6.} Theorem} The $\BdR$-linear extension of $\rho$ is a filtered quasi-isomorphism: for any $X\in\CV ar_{\bar{K}}$ one has $$\rho : R\Gamma_{\!\dR}(X)\otimes_{\bar{K}} \BdR \iso R\Gamma_{\!\et}  (X,\Bbb Q_p )\otimes_{\Bbb Q_p} \BdR . \tag 3.6.1   $$\endproclaim

\demo{Proof} (a) {\it The case of $X=\Bbb G_{m}=\Bbb G_{m\bar{K}}$}:  The $\bar{K}$-line 
 $H^1_{\dR} (\Bbb G_{m} )=\gr^1_F H^1_{\dR} (\Bbb G_{m} )$ is generated by $d\log t$. The $\Bbb Z_p$-line $H^1_{\et} (\Bbb G_{m},\Bbb Z_p )(1)= H^1_{\et} (\Bbb G_{m},\Bbb Z_p (1))$ is generated by the class $cl( \fk )$ of  the Kummer $\Bbb Z_p (1)$-torsor $\fk =\limleft\fk_n$, $\fk_n :=(t^{1/p^n})$. Due to the canonical identification $ \Bbb C_p (1)\iso \fm_{\dR}/\fm_{\dR}^2 =
 \gr^1_F \BdR $, see 1.4, 1.5, we can view  $cl( \fk )$ as a generator of the $\Bbb C_p$-line 
$H^1_{\et} (\Bbb G_{m},\Bbb Q_p )\otimes  \gr^1_F \BdR$.

\proclaim { Lemma} One has $\gr^1_F (\rho )(d\log t )= cl( \fk )\in H^1_{\et}  (\Bbb G_{m},\Bbb Q_p )\otimes\gr^1_F \BdR $.\endproclaim

\demo{Proof of Lemma} We make a mod $p^n$ computation. 
Consider the ss-pair $(\Bbb G_{m },\bar{\Bbb G}_{m})$,    $ \bar{\Bbb G}_{m} :=   \Bbb P^1_{O_{\bar{K}}}$. One has $\gr^1_F L\Omega^\cdot\!\hat{_{(\Bbb G_m }}{}_{,\bar{\Bbb G}_m)}= 
  \Omega^1_{(\Bbb G_{m },\bar{\Bbb G}_{m})}[-1]$, so 
  $ d\log t \in  \Gamma (  \bar{\Bbb G}_{m},    \Omega^1_{(\Bbb G_{m },\bar{\Bbb G}_{m})})$ is a 1-cocycle in $\gr^1_F R\Gamma^\natural_{\!\dR}(\Bbb G_m )$.  As in 1.1, set $C_n :=\CC one \, (p^n : \Bbb Z \to \Bbb Z )$.  Let   $d\log cl(\fk_n )$ be  the image of the class
 $cl(\fk_n )$  of $\fk_n$ by the composition 
 $H^1(\Bbb G_{m\,\et} ,\mu_{p^n})\to H^1(\Bbb G_{m\,\et},\gr_F^1 \AdR \otimes C_n )\to \gr^1_F  R\Gamma^\natural_{\!\dR} (\Bbb G_{m})\otimes C_n $, where the first arrow comes from the coefficient maps $\mu_{p^n} \buildrel{d\log}\over\lra \Omega_{O_{\bar{K}}\, p^n} \hra
 \Omega_{O_{\bar{K}}}[-1] \otimes C_n = \gr_F^1 \AdR \otimes C_n$. 
To prove the lemma, we will show that the image of $d\log t$ by the embedding  $ \gr^1_F  R\Gamma^\natural_{\!\dR} (\Bbb G_m )\hra  \gr^1_F  R\Gamma^\natural_{\!\dR} (\Bbb G_m )\otimes C_n$ is homologous to   $d\log cl(\fk_n )$.

Let $\Bbb G\tilde{_m}$ be a copy of $\Bbb G_m$ with parameter $\tilde{t}$, and $\pi: \Bbb G\tilde{_m}\to \Bbb G_m$ be the projection $t=\tilde{t}^{p^n}$. Thus $\Bbb G\tilde{_m}/\Bbb G_m$ is our  $\mu_{p^n}$-torsor  $\fk_n$, so  $cl(\fk_n )$ is represented by a \v Cech $\mu_{p^n}$-cocycle $c (\fk_n )$ for the \'etale covering $ \Bbb G\tilde{_m}/\Bbb G_m$. The corresponding \v Cech hypercovering is the twist of $ \Bbb G\tilde{_m}$ by the universal $\mu_{p^n}$-torsor $\ft_n$ over the classifying simplicial space $B_{\mu_{p^n}\cdot}$, so
 for any sheaf $\CF$ the \v Cech complex of $ \Bbb G\tilde{_m}/\Bbb G_m$
  with coefficients in $\CF$ is the cochain complex $C^\cdot (\mu_{p^n}, \Gamma ( \Bbb G\tilde{_m},\CF ))$ for $\mu_{p^n}$ acting on sections by the translations. The 1-cocycle $c (\fk_n )$ is the identity map $\mu_{p^n}\to \mu_{p^n}=\Gamma ( \Bbb G\tilde{_m},\mu_{p^n})$.

Our $\pi$ extends to the h-covering of semi-stable pairs $ ( \Bbb G\tilde{_m}, \bar{ \Bbb G}\tilde{_m}   )\to ( \Bbb G_m ,  \bar{\Bbb G}_m )$, and the \v Cech hypercovering extends 
to a hypercovering in $\CV ar^{\ss}_{\bar{K}}$ which is
the $\ft_n$-twist of $ ( \Bbb G\tilde{_m}, \bar{ \Bbb G}\tilde{_m})$. So one has a canonical map $ C^\cdot (\mu_{p^n}, \Gamma (\bar{ \Bbb G}\tilde{_m}, \Omega^1_{( \Bbb G\tilde{_m}, \bar{ \Bbb G}\tilde{_m})}))[-1] \to \gr^1_F  R\Gamma^\natural_{\!\dR} ( \Bbb G_{m}  )$, hence $  C^\cdot (\mu_{p^n}, \Gamma (\bar{ \Bbb G}\tilde{_m}, \Omega^1_{( \Bbb G\tilde{_m}, \bar{ \Bbb G}\tilde{_m})}))[-1]\otimes C_n  \to 
 \gr^1_F  R\Gamma^\natural_{\!\dR} ( \Bbb G_{m}  ) \otimes C_n $. 
Both $d\log t$ and $d\log c (\fk_n )$ are  1-cocycles in  $  C^\cdot (\mu_{p^n}, \Gamma (\bar{ \Bbb G}\tilde{_m}, \Omega^1_{( \Bbb G\tilde{_m}, \bar{ \Bbb G}\tilde{_m})}))[-1]\otimes C_n$:
namely, $d\log t \in C^0 (\mu_{p^n},  \Gamma (\bar{ \Bbb G}\tilde{_m}, \Omega^1_{( \Bbb G\tilde{_m}, \bar{ \Bbb G}\tilde{_m})})[-1])$ and
$d\log c(\fk_n )\in C^1 (\mu_{p^n}, \Omega_{O_{\bar{K}}\,p^n} )\subset C^1 (\mu_{p^n}, \Gamma (\bar{ \Bbb G}\tilde{_m}, \Omega^1_{( \Bbb G\tilde{_m}, \bar{ \Bbb G}\tilde{_m})})[-1]\otimes C_n )$.
Their difference is the differential of the 0-cochain $d\log \tilde{t} \in C^0 (\mu_{p^n},  \Gamma (\bar{ \Bbb G}\tilde{_m}, \Omega^1_{( \Bbb G\tilde{_m}, \bar{ \Bbb G}\tilde{_m})}))
\subset C^0 (\mu_{p^n},  \Gamma (\bar{ \Bbb G}\tilde{_m}, \Omega^1_{( \Bbb G\tilde{_m}, \bar{ \Bbb G}\tilde{_m})})[-1]\otimes C_n )$,  q.e.d.
\qed  \enddemo 
 
We see that  for $X=\Bbb G_m$ the map $\rho$ of (3.6.1) is a filtered quasi-isomorphism. It provides a canonical generator $\rho (d\log t)/cl (\fk )$ of $\fm_{\dR}(-1)$. Thus we have a canonical identification between $H^\cdot_{\et}  (X,\Bbb Q_p (n))\otimes \BdR$ and $H^\cdot_{\et}  (X,\Bbb Q_p )\otimes \BdR$ with the filtration shifted by $n$, etc.

\remark{Remark} The above generator is equal to  Fontaine's generator, see \cite{F3} 1.5.4. Indeed, they coincide modulo $\fm_{\dR}^2 (-1)$ by the lemma, and both are Gal$(\bar{K}/K)$- invariant (see Remark in 3.5). Since $H^0 (\text{Gal}(\bar{K}/K) , \fm_{\dR}^{2}(-1))=0$, we are done. 
\endremark
\enspace

(b) {\it Compatibility with the Gysin maps}: Let $i: Y\hra X$ be a closed 
codimension 1 embedding of smooth varieties. It yields the Gysin isomorphisms $i_{*\dR} : R\Gamma_{\!\dR }(Y) \iso  R\Gamma_{\!\dR\, Y}(X)(1)[2]:=  \CC one ( R\Gamma_{\!\dR }(X)\to \Gamma_{\!\dR }(X\smallsetminus Y))(1)[1]$,  
$i_{*\Bbb Q_p} : R\Gamma_{\!\et}  (Y,\Bbb Q_p ) \iso  R\Gamma_{\!{\et}\, Y}(X,\Bbb Q_p )(1)[2]$. Let us show that $\rho$ commutes with the Gysin maps. \newline

 Consider the deformation to the normal cone diagram $$\spreadmatrixlines{1\jot}
\matrix
\CL
 &\hra&
 X\tilde{_{\Bbb A^1}} & \hookleftarrow &X   \\  
\uparrow  &&\uparrow&& \uparrow  \\
Y &  \hra & Y_{\Bbb A^1} &\hookleftarrow & Y. 
    \endmatrix
\tag 3.6.2$$
Here $Y_{\Bbb A^1}= Y\times \Bbb A^1$, $X\tilde{_{\Bbb A^1}}$ is $X\times \Bbb A^1$ with $Y\times\{0\}$ blown up, the left arrow is the zero section of the normal bundle $\CL$ to $Y$ in $X$, and the bottom embeddings are $y\mapsto (y,0), (y,1)$. It yields a commutative diagram of the de Rham cohomology
$$\spreadmatrixlines{1\jot}
\matrix
R\Gamma_{\!\dR\, Y }( \CL )(1)[2]
 &\leftarrow&R\Gamma_{\!\dR\, Y_{\Bbb A^1} }
 (X\tilde{_{\Bbb A^1}})(1)[2] & \to &   R\Gamma_{\!\dR\, Y}(X)(1)[2]   \\  
\uparrow  &&\uparrow&& \uparrow  \\
R\Gamma_{\!\dR}(Y) &  \leftarrow & R\Gamma_{\!\dR}(Y_{\Bbb A^1}) &\to  &  R\Gamma_{\!\dR}(Y),   \endmatrix
\tag 3.6.3$$ 
 where the vertical arrows are the Gysin isomorphisms and the horizontal ones are pullbacks. There is 
 a similar diagram for the $\Bbb Q_p$-cohomology. The horizontal maps are filtered quasi-isomorphisms, so, since $\rho$ is compatible with pullbacks, we see that
the Gysin compatibility for $Y\hra X$ amounts to one for $Y\hra \CL$. 

So we 
 can assume that $X$ is a line bundle $\CL$ over $Y$ and $i$ its zero section. 
  Now the source of both $i_*$'s are dg algebras, the targets are modules over them (due to the projection $\CL\to Y$), and $i_*$'s are morphisms of modules.
Thus it suffices to check that $\rho$ identifies the images of 1. The assertion is local with respect to $Y$, hence we can assume that $\CL$ is trivial. By base change, we reduced to the case when $Y$ is a point, where we are done by (a).

(c) {\it The case of a smooth projective $X$}: Let us check that the morphism of bigraded rings $\gr_F^\cdot \rho^* : \gr_F^\cdot H^*_{\dR}(X)\otimes_{\bar{K}}  \gr_F^\cdot \BdR \to H^*_{\et} (X,\Bbb Q_p )\otimes_{\Bbb Q_p}\gr_F^\cdot \BdR$ is an isomorphism.  
It is an isomorphism for $*=0$. By (b), $\gr_F^1 \rho^2$ identifies the classes $c$ of a hyperplane section. Since the product with $c^{\dim X}$ identifies $H^0$ and $H^{2\dim X}$,   $\gr_F^\cdot \rho^{2\dim X}$  is an isomorphism. Therefore, since $\gr_F^\cdot \rho^*$ is compatible with the Poincar\'e pairing for classes of opposite degrees and the latter is nondegenerate, $\gr_F^\cdot \rho^*$ is injective. Since  $\dim_{\bar{K}}H^*_{\dR}(X)=\dim_{\Bbb Q_p} H^*_{\et}  (X,\Bbb Q_p )$, we are done.

(d) {\it The case of $X=\bar{X}\smallsetminus D$, where $\bar{X}$ is smooth projective, $D$ is a   strict normal crossings divisor}:
Let $Y$ be an irreducible component of $D$, $D'$ be the union of the other components; set $X':=\bar{X}\smallsetminus D'$, $Y':=Y\smallsetminus D'$. By induction by the number of components $D$ (starting with (c)), we can assume that the theorem holds  for $X'$ and $Y'$. By (b), $\rho$ provides a morphism between the exact Gysin triangles for $(Y', X')$. It is a filtered quasi-isomorphism on the $X'$ and $Y'$ terms; hence it is a filtered quasi-isomorphism on the $X$ term, q.e.d.

(e) {\it The case of arbitrary $X$}: If $Y_\cdot /X$ is any h-hypercovering of $X$, then the canonical map $R\Gamma_{\!\dR}(X)\to R\Gamma_{\!\dR}(Y_\cdot )$ (which is the total  complex of the cosimplicial system of filtered complexes $R\Gamma_{\!\dR}(Y_i )$) is a filtered quasi-isomorphism by the construction of $R\Gamma_{\!\dR}$, and $R\Gamma_{\!\et}(X,\Bbb Q_p )\iso R\Gamma_{\!\et}(Y_\cdot ,\Bbb Q_p )$ by cohomological descent (see Remark in 2.4). Thus
if $\rho$ is a filtered quasi-isomorphism for every $Y_i$, then it is  a filtered quasi-isomorphism for $X$. We are done, since, by de Jong (or Hironaka), one can find $Y_\cdot /X$ with $Y_i$ as in (d). 
\quad\qed  \enddemo 

\remark{Remark} $\rho$ is compatible with Chern classes of vector bundles: Indeed,  $c_i (E)$ are determined in the usual way by $c_1 (\CO (1)_{\Bbb P (E)})$, so it suffices to show that $\rho$ identifies $c_1$'s of line bundles.  Notice that
the construction of $\rho$ extends tautologically to simplicial schemes. By (a) above, $\rho$ identifies the de Rham and \'etale Chern classes of the universal line bundle over the classifying simplicial scheme $B_{\Bbb G_m  \cdot}$. For  a line bundle $\CL$ on $X$, choose a finite open covering $\{ U_i \}$ of $X$ such that $\CL$ is trivial on $U_i$; let $\pi : X\tilde{_\cdot}\to X$ be the  \v Cech hypercovering. Since  $\pi$ yields an isomorphism between the cohomology,  it suffices to check that $\rho$ identifies the Chern classes of $\pi^*\CL$. This is true since $\pi^* \CL$ is the pullback of the universal line bundle by a map $X\tilde{_\cdot} \to B_{\Bbb G_m  \cdot}$. 
\endremark

\head 4. Proof of the Poincar\'e Lemma. \endhead

\subhead{\rm 4.1}\endsubhead Pick any $(V,\bar{V})\in\CV ar^{\ss}_{\bar{K}}$.

\proclaim {Proposition} One has  $\text{L}_{(V,\bar{V})}\iso \Omega_{(V,\bar{V})}$, the $\CO_{\bar{V}}$-module  $\Omega_{\langle V,\bar{V}\rangle }:= \Omega_{(V,\bar{V})/ (\bar{K},O_{\bar{K}})}$ is locally free of finite rank, and there is a canonical short exact sequence   $$
0\to\CO_{\bar{V}}\otimes_{O_{\bar{K}}} \Omega_{O_{\bar{K}}}\to\Omega_{(V,\bar{V})}\to \Omega_{\langle V,\bar{V}\rangle }
\to 0. \tag 4.1.1$$ 
\endproclaim

\demo{Proof} We can assume that $V$ is connected, so  $(V,\bar{V})$ is the base change of a semi-stable $K$-pair $(U,\bar{U})$ as in 2.2(c), i.e., $ (V,\bar{V})=(U_{\bar{K}}, \bar{U}_{O_{\bar{K}}})$. For any
 finite extension $K'$ of $K_U$, consider an arithmetic $K$-pair $(U_{K'},\bar{U}_{O_{K'}}):= (U\otimes_{K_U}\! K', \bar{U}\otimes_{O_{K_U}}\! O_{K'})$. Set $\Omega_{\langle U,\bar{U}\rangle}:= \Omega_{(U,\bar{U})/(K_U,O_{K_U})}$, $\Omega_{\langle U_{K'},\bar{U}_{O_{K'}}\rangle}:= \Omega_{(U_{K'},\bar{U}_{O_{K'}})/(K',O_{K'})}$.

\proclaim{Lemma} The log scheme $(U_{K'},\bar{U}_{O_{K'}})$ coincides with  the  
pullback  of  $(U,\bar{U})$ by the map $\Spec\, (K',O_{K'} )\to \Spec\, (K_U,O_{K_U} )$ in the category of log schemes. \endproclaim

Assume the lemma for a moment. The map $(U,\bar{U})\to\Spec\, (O_{K_U} ,K_U)$ is log smooth and integral; by the lemma,  $ (U_{K'},\bar{U}_{O_{K'}})\to \Spec\, (O_{K'}, K')$ enjoys the same properties.  So, by \cite{Ol} 8.34, L$_{ (U_{K'},\bar{U}_{O_{K'}})/ (K',O_{K'})}\iso \Omega_{\langle U_{K'},\bar{U}_{O_{K'}}\rangle}=O_{K'}\otimes_{O_{K_U}}\Omega_{\langle U,\bar{U}\rangle}$, which is a locally free $\CO_{ \bar{U}_{O_{K'}}}$-module of finite rank. Since $\text{L}_{(K',O_{K'})}
\iso \Omega_{(K',O_{K'})}$ (see the proof of Lemma in 3.2), 
 the canonical exact triangle  (\cite{Ol} 8.18) 
 $\CO_{\bar{U}_{O_{K'}}}\otimes_{O_{K'}}\text{L}_{(K',O_{K'})}\to \text{L}_{(U_{K'},\bar{U}_{O_{K'}})}\to \text{L}_{ (U_{K'},\bar{U}_{O_{K'}})/ (K',O_{K'})}$ reduces to the short exact sequence $0 \to \CO_{\bar{U}_{O_{K'}}}\otimes_{O_{K'}}\Omega_{(K',O_{K'})} \to \Omega_{(U_{K'},\bar{U}_{O_{K'}})}\to \Omega_{\langle U_{K'},\bar{U}_{O_{K'}}\rangle}
\to 0$.  Pass to the limit by all $K'\subset \bar{K}$ and  use the lemma in 3.2; we are done.  \qed
  \enddemo 

\demo{Proof of Lemma} The underlying scheme of the pullback log scheme is $\bar{U}_{O_{K'}}$. Let us show that its log structure map $\CM\to \CO_{\bar{U}_{O_{K'}}}\cap   j_* \CO^\times_{U_{K'}}$  is an isomorphism. The assertion is \'etale local, so
we can assume that $\bar{U}$ is \'etale over $\Spec\, O_{K_U} [ t_a ,  t_b , t_c ]/$ $(\Pi t_a -\pi_{K_U} )$, where $a$, $b$, $c$ are in finite sets $A$, $B$, $C$, $\pi_{K_U}$ is a uniformizing parameter in $O_{K_U}$, and $U$ is the subscheme where all $t_a$,  $t_b$ are invertible. The log structure of
 $(U,\bar{U})$ is fine with a chart $\Bbb N [A\sqcup B] \to \CO_{\bar{U}}$, which sends generators   $m_a, m_b$ of $\Bbb N[A\sqcup B]$ to
   $ t_a, t_b$. Therefore\footnote{We replace one $t_a \in\CO( \bar{U}_{O_{K'}})$ by $t_a \pi_{K'}^e /\pi_{K_U}$. } $\bar{U}_{O_{K'}}$ is \'etale over $\Spec\, O_{K'} [ t_a ,  t_b ,  t_c ]/(\Pi t_a -\pi^e_{K'} )$, where $e$ is the ramification index of $K'/K_U$,  $\pi_{K'}$ is a uniformizing parameter in $O_{K'}$,  and the log structure $\CM$ has a chart $M_{A,B}\to \CO_{\bar{U}_{O_{K'}}}$, where $M_{A,B}$ is the quotient of    $\Bbb N [A\sqcup B]\oplus \Bbb N$ modulo the relation $\Pi m_a =m^e_\pi$ ($m_\pi$ is  the generator of the last summand $\Bbb N$),  the chart is  $m_a, m_b, m_\pi \mapsto t_a, t_b, \pi_{K'}$. Consider an embedding
 $M_{A,B}\hra M^w_{A,B}:= e^{-1}\Bbb N [A] \oplus \Bbb N [B]$, $m_a, m_b, m_\pi \mapsto m_a ,m_b, \Pi m_a^{1/e}$. Its image is formed by those $\Pi m_a^{n_a /e} \Pi m_b^{n_b}$, $n_a ,n_b \in \Bbb N$, such that $n_a -n_{a'}\in e\Bbb Z$ for any $a,a'\in A$. Thus $M_{A,B}$ is saturated.
Now  the log scheme $(\bar{U}_{O_{K'}},\CM )$ is evidently log regular in the sense of \cite{K2} 2.1, hence $\CM \iso \CO_{\bar{U}_{O_{K'}}}\cap j_* \CO^\times_{U_{K'}}$ by \cite {K2} 11.6, q.e.d.\footnote{I am grateful to Luc Illusie for the proof.}

The reference to \cite{K2} can be replaced by the next
 explicit argument: It suffices to show that the map of sheaves $\CM/\CO^\times_{\bar{U}_{O_{K'}}} \to ( \CO_{\bar{U}_{O_{K'}}  }\cap j_* \CO^\times_{U_{K'}})/\CO_{\bar{U}_{O_{K'}} }^\times$ is an isomorphism. The r.h.s.~is the sheaf $\CD$ of effective Cartier divisors supported on $\bar{U}_{O_{K'}}   \smallsetminus U_{K'}$. Let $\CD^w \supset \CD$ be the sheaf of the corresponding effective Weil divisors. For $x\in \bar{U}_{O_{K'}} $, the fiber 
$(\CM/\CO^\times_{\bar{U}_{O_{K'}}})_x$  is the quotient $M_{A_x,B_x}$ of $M_{A,B}$, where $A_x \subset A$, $B_x \subset B$ consist of those $a,b$ such that $t_a ,t_b$ vanish at $x$. The map $M_{A_x,B_x}\to \CD_x$ extends to an
isomorphism $M^w_{A_x,B_x}\iso \CD^w_x$, which identifies a generator $m^{1/e}_a $
with  the {\it reduced} divisor $D_a$ of $t_a$, $m_b$ with $D_b :=\text{div} ( t_b )$. Thus $M_{A_x,B_x}\hra \CD_x$. To show that $\hra$ is an isomorphism, we need to check that if $D=\Sigma n_a D_a +\Sigma n_b D_b$ is a Cartier divisor at $x$, then $n_a -n_{a'}\in e\Bbb Z$ for any $a,a'\in A_x$. We can assume that $A=\{ a,a'\}$, $B=C=\emptyset$, so $\bar{U}_{O_{K'}}  $ is a semi-stable curve over $O_{K'}$. The  exceptional divisor of its minimal desingularization $\tilde{U}_{O_{K'}} $  is a chain of $e-1$ projective lines $P_1 , \ldots , P_{e-1}$ with self-intersection indices $(P_i ,P_i )= -2$.  
Let $\tilde{D}= n_a D^\flat_a + n_1 P_1 + \ldots + n_{e-1}P_{e-1} + n_{a'}D^\flat_{a'}$ be the pullback of $D$ to $\tilde{U}_{O_{K'}} $; here $D_a^\flat$, $D_{a'}^\flat$ are strict transforms of $D_a$, $D_{a'}$. One has
 $(\tilde{D},P_i )=0$, i.e., $n_{i-1}-2n_i + n_{i+1}=0$ or $n_i - n_{i-1}=n_{i+1}-n_i$, where $n_0 := n_a$, $n_e := n_{a'}$. Thus $n_{a'}-n_a =e(n_1 - n_a)\in e\Bbb Z$, and we are done.
 \qed
 \enddemo

\subhead{\rm 4.2}\endsubhead  Set
      $ \Omega^a_{\langle V,\bar{V}\rangle }:= \Lambda^a_{\CO_{\bar{V}}}  \Omega_{\langle V,\bar{V}\rangle }$. Consider (4.1.1) as
a  2-step filtration  on $\Omega_{(V,\bar{V})}$; it  splits locally since $\Omega_{\langle V,\bar{V}\rangle}$ is locally free. Passing to derived exterior powers, we  get for any $m$ a finite increasing filtration $I_\cdot$  on 
 $\gr^m_F L\Omega^\cdot \!\hat{_{(V,}}{}_{\bar{V})}=   (L\Lambda^m_{\CO_{\bar{V}}} \Omega_{(V,\bar{V})})[-m]$ with $\gr^I_a \gr^m_F 
 L\Omega^\cdot \!\hat{_{(V,}}{}_{\bar{V})}
 = \Omega^{a}_{\langle V,\bar{V}\rangle} \otimes_{O_{\bar{K}}} \gr^{m-a}_F \AdR [-a]$, hence 
 on $\gr^m_F R\Gamma^\natural_{\!\dR} (V,\bar{V})$ with 
 
 $$\gr^I_a  \gr^m_F R\Gamma^\natural_{\!\dR} (V,\bar{V})
 = R\Gamma (\bar{V},  
 \Omega^{a}_{\langle V,\bar{V}\rangle}) \otimes^L_{O_{\bar{K}}} \gr^{m-a}_F \AdR[-a]. \tag 4.2.1$$ 
 
Let $\CG^a$ be the h-sheafification (see (2.6.1)) of the complex of presheaves  $(V,\bar{V})\mapsto R\Gamma (\bar{V}, \Omega^{a}_{\langle V,\bar{V}\rangle})$ on $\CV ar^{\ss}_{\bar{K}}$. This is a complex of h-sheaves of $O_{\bar{K}}$-modules on $\CV ar_{\bar{K}}$. Its cohomology $H^b \CG^a$ is h-sheafification of the presheaf $(V,\bar{V})\mapsto H^b (\bar{V}, \Omega^{a}_{\langle V,\bar{V}\rangle})$ on $\CV ar^{\ss}_{\bar{K}}$.
Our
 $I_\cdot$ is a filtration on the presheaf $(V,\bar{V})\mapsto \gr^m_F R\Gamma^\natural_{\!\dR}(V,\bar{V})$; passing to h-sheafification, we get a finite filtration
$I_\cdot$ on $\gr^m_F \CA^\natural_{\dR}$ with span $[0,m]$ and $$
  \gr^I_{a}\,\gr^m_F \CA^\natural_{\dR} = \CG^{a}\otimes^L_{O_{\bar{K}}} \gr^{m-a}_F \AdR [-a]. \tag 4.2.2$$ Notice that the bottom cohomology $H^0$ of the bottom term $I_0 =\gr^I_0$ is the constant sheaf $O_{\bar{K}}$ and  $\CC one (\gr^m_F \AdR\to \gr^m_F \CA^\natural_{\dR} )=\gr^m_F \CA^\natural_{\dR} /H^0 I_{0}$. Therefore, by (4.2.2), the Poincar\'e lemma follows from the next assertion:

 \proclaim{Theorem} The cohomology $H^b \CG^a$ are h-sheaves
 of $\Bbb Q$- (hence $\bar{K}$-) vector spaces for $(a,b)\neq (0,0)$.
  \endproclaim
 
 \remark{Remark} 
 The p-divisiblity of $H^{b}\CG^0$, $b\neq 0$, was first proved  by Bhatt  \cite{Bh1} 8.0.1. \endremark
\remark
{ Exercise} Consider a presheaf $(V,\bar{V})\mapsto R\Gamma (\bar{V}, 
 L\Omega^\cdot \!\hat{_{(V,}}{}_{\bar{V})/ (\bar{K},O_{\bar{K}})})$; let $\CA_{\dR}^{naive}$ be its  h-sheafification. One has an evident map $\CC one (F^1 \AdR\to\CA^\natural_{\dR} ) \to\CA_{\dR}^{naive}$. Show that
 the theorem implies that it is a filtered quasi-isomorphism, i.e., the
  triangle $F^1 \AdR\to\CA^\natural_{\dR} \to\CA_{\dR}^{naive}$ is exact in the filtered derived category of h-sheaves.
\endremark
  
\subhead{\rm4.3}\endsubhead We deduce the above theorem from a more concrete assertion. As in 4.1, for an ss-pair $(U,\bar{U})$ over $K$ we have the
locally free $\CO_{\bar{U}}$-module of log differentials $\Omega_{\langle U,\bar{U}\rangle}:=\Omega_{(U,\bar{U})/(K_U ,O_{K_U} )}$ and its exterior powers
$\Omega^a_{\langle U,\bar{U}\rangle}:=\Lambda^a \Omega_{\langle U,\bar{U}\rangle}$.

Let $f: (U',\bar{U}')\to (U,\bar{U})$ be a map in $\CV ar^{\ss}_K$ or $\CV ar^{\ss}_{\bar{K}}$.
We  say that $f$ is {\it (Hodge) $p$-negligible} if  the morphisms 
$(\tau_{>0} R\Gamma 
(\bar{U}, \CO_{\bar{U}}))\otimes^L \Bbb Z/p \to (\tau_{>0} R\Gamma (\bar{U}', \CO_{\bar{U}'}))\otimes^L \Bbb Z/p$ and
$R\Gamma 
(\bar{U}, \Omega^a_{\langle U,\bar{U}\rangle})\otimes^L \Bbb Z/p \to R\Gamma (\bar{U}', \Omega^a_{\langle U',\bar{U}'\rangle})\otimes^L \Bbb Z/p$, $a>0$, in $D^b (O_{K_U} /p )$, resp.~$D^b (O_{\bar{K}} /p )$,  vanish. 

\remark{Remark} For $(U,\bar{U})\in \CV ar^{\ss}_K$ and a point $K_U \to\bar{K}$, one has 
$R\Gamma (\bar{U}_{O_{\bar{K}}}, \Omega^a_{\langle U_{\bar{K}},\bar{U}_{O_{\bar{K}}}\rangle})=R\Gamma 
(\bar{U}, \Omega^a_{\langle U,\bar{U}\rangle})\otimes^L_{O_K} O_{\bar{K}}$.
Therefore the base change functor $\CV ar^{\ss}_K \to\CV ar^{\ss}_{\bar{K}}$  (see (2.2.2)) preserves $p$-negligible  maps.\endremark

\proclaim{Theorem} Every  $U\in\CV ar^{\ss}_K$ 
 admits a p-negligible h-covering. Ditto for $\bar{K}$-pairs.\endproclaim

The theorem implies the one in 4.2: Indeed, the  $\bar{K}$-assertion shows that one has
 $(\tau_{>0}\CG^0 )\otimes^L \Bbb Z/p =0$ and $\CG^a \otimes^L \Bbb Z/p =0$ for $a>0$; since for a complex $\CG$ the multiplication by $p$ on $H^\cdot \CG$ is invertible if and only if $\CG\otimes^L \Bbb Z/p =0$, we are done. Thus it yields the Poincar\'e lemma. 
 
The above remark shows that the $K$-version of the theorem implies the $\bar{K}$-one. The proof of the $K$-version  takes the rest of the section.

\subhead{\rm 4.4}\endsubhead
{\it For the rest of \S 4, ``pair" means ``arithmetic $K$-pair" (see 2.2)}. We need further input from de Jong. 
 A morphism $f: (C,\bar{C}) \to (S,\bar{S})$  of pairs is said to be a {\it family of pointed curves} (over $(S,\bar{S})$) if the map $\bar{C}_S := f^{-1}(S)\to S$ is smooth  of relative dimension 1 with irreducible geometric fibers, and
$D_{fS}:=\bar{C}_S \smallsetminus C$, viewed as a reduced scheme, is \'etale over $S$. Such an  $f$ is {\it semi-stable} if, in addition,  
$ \bar{C}/\bar{S}$ is a semi-stable family of curves, and the closure $D_f$ of $D_{fS}$ in $\bar{C}$ (the {\it  horizontal divisor}), viewed as a reduced scheme, is \'etale over $\bar{S}$ and  intersects each fiber of $f$ at smooth points. A  section $e: (S,\bar{S})\to (C,\bar{C})$  of $f$ is said to be {\it nice} if $e(\bar{S})$ intersects fibers of $f$ at smooth points and $D_f \cap e(\bar{S})=\emptyset$. Families of pointed curves over $(S,\bar{S})$ form a category $\CC_{(S,\bar{S})}$ in the obvious manner, and a morphism of bases $ \psi : (S',\bar{S}')\to (S,\bar{S})$ yields an evident pullback functor
$ \CC_{(S,\bar{S})}\to\CC_{(S',\bar{S}')}$ which preserves semi-stable families. A morphism $ f'\to f$ in $\CC_{(S,\bar{S})}$ is called {\it alteration} if $(C',\bar{C}')$ is an alteration of $(C,\bar{C})$; it is
 a {\it semi-stable alteration} (of $f$) if, in addition, $f'$ is semi-stable. 

\proclaim {Theorem}   (a) Any family $f: (C,\bar{C}) \to (S,\bar{S})$  of pointed curves with $f: \bar{C}\to\bar{S}$  projective admits  
a semi-stable alteration $f'$
h-locally over $(S,\bar{S})$.
\newline
(b) One can find $f'$ as above which has a nice section $e$. 
 Moreover, for a given  closed subscheme $P\subset \bar{C}$ such that $f(P)=\bar{S}$ and $P\cap \bar{C}_S \subset C$,  one can find $e$ such that the map $\bar{C}'\to \bar{C}$ sends $e(\bar{S})$ to $P$. 
 \newline 
(c) For any semi-stable family of pointed curves $f: (C,\bar{C}) \to (S,\bar{S})$
with $(S,\bar{S})$ a strict ss-pair, there exists a semi-stable alteration   $m: (C,\tilde{C})\to (C,\bar{C})$ of $f$ with $m|_C =\id_C$  such that  $m :\tilde{C}\to \bar{C}$ is an isomorphism over smooth points of $f$ and $  (C,\tilde{C})$ is an ss-pair.\endproclaim

\demo{Proof}   (c) is \cite{dJ1} 3.6. 
(a) follows from  \cite{dJ2} 2.4 (i),(ii) except that de Jong does not care to control the domain of smoothness of the semi-stable alteration of $f$.  A miniscule modification of his argument permits to do this; see Appendix 1. Alternatively, (a) follows directly from a far more precise result of Temkin \cite{T} 1.5.\footnote{To bring our datum to Temkin's setting, one flattens $\bar{C}/\bar{S}$ and $D_f /\bar{S}$ using \cite{RG} 5.2.2 and replaces $\bar{S}$ by its normalization.}

Let us check (b). Every pair has a canonical alteration by the union of normalizations of its irreducible components, so we assume all the way that $\bar{S}$ is normal and irreducible. Since $P$ as in (b)  exists h-locally on $(S,\bar{S})$,\footnote{Indeed, one can cover $S$ by  open subsets  $\{ S_\alpha \}$ such that $P$ exists for $(C_{S_\alpha},\bar{C})\to (S_\alpha ,\bar{S})$.} we can assume  it is given. Replacing  $(S,\bar{S})$ by its alteration $(P_S ,P)$, we get a section $e$ of $f$  with image in $P$. Set $C^{\flat}:= C\smallsetminus e(S)$. Then $(C^\flat ,\bar{C})\to (S,\bar{S})$ is a family of pointed curves; let $f^{\prime\flat} : (C^{\prime \flat},\bar{C}')\to  (S,\bar{S})$ be its semi-stable alteration as in (a). Let $D_e$ be the closure in $ \bar{C}'$ of the preimage $D_{eS}$ of $e(S)$. Then $D_e$ is an \'etale covering of $\bar{S}$.\footnote{Since $D_{f^{\prime\flat}}$ is  \'etale over $\bar{S}$, $D_{eS}$ is  open and closed in $D_{f^{\prime\flat}S}$, $f^{\prime\flat} (D_{eS})=S$,
and $\bar{S}$ is normal.} Let $C' $ be the preimage of $C\subset \bar{C}$ in $\bar{C}'$; then $(C',\bar{C}')\to (S,\bar{S})$ is a semi-stable  alteration of $(C,\bar{C})\to  (S,\bar{S})$. Replacing $(S,\bar{S})$ by its alteration $(D_{eS}, D_e )$, we get a nice section of $(C',\bar{C}')$ which sits over $e$, hence over $P$.
 \quad\qed  \enddemo 
 
\remark{Remark} In  (c), every  nice section of $(C,\bar{C})$ lifts to a nice section of $(C,\tilde{C})$.\endremark

\proclaim {Corollary}  Any pair $(U,\bar{U})$ has an h-covering by ss-pairs $(C,\bar{C})$, $\dim C=\dim U$, for which there is a semi-stable family of pointed curves $f : (C ,\bar{C} )\to (S ,\bar{S} )$ with a nice section such that $(S,\bar{S} )$ is  an ss-pair and $C$ is  affine over $S$ (i.e.,  $f(D_{f})= \bar{S}$).\endproclaim

\demo{Proof} It suffices to find an h-covering of $(U,\bar{U})$ by pairs $(C,\bar{C})$ with $\dim C=\dim U$ for which there exists a family of pointed curves $f : (C ,\bar{C} )\to (S ,\bar{S} )$ with $C$ affine over $S$ and projective $\bar{S}$, $\bar{C}$. The theorem transforms it then, with an input from Remark (i) in 2.5 to preserve the dimension 
and de Jong's theorem in 2.3 to alter $(S,\bar{S})$ from (b) into a strict ss-pair, into a datum with all  promised properties. 

By de Jong's theorem in 2.3, we can assume that 
$(U,\bar{U})$ is an ss-pair and $\bar{U}$ is projective and irreducible;\footnote{We  only need that $\bar{U}$ is projective and normal, and that $U$ is smooth.}  set $d=\dim U$. 
Pick any closed point $u\in U$. It suffices to find an open neighborhood $U'\subset U$ of $U$, an alteration $(C,\bar{C})$ of $(U',\bar{U})$, and a family of curves $f: (C,\bar{C})\to (S ,\bar{S} )$ such that $f(D_{f})=\bar{S}$.

Embed $\bar{U}$ into a projective space $\Bbb P^N_{O_K}$. 
By Bertini, there is a plane $H\subset \Bbb P^N$ defined over $K$ of codimension $d$ such that $u\notin H$, $H$ intersects $\bar{U}_K$ transversally, $H\cap \bar{U}_K \subset U$,
and the  codimension $d-1$ plane which contains $H$ and $u$, is transversal to $\bar{U}_K$ and $\bar{U}_K \smallsetminus U$. Let $m: \bar{C}\to \bar{U}$ be the blowup at $\bar{U}\cap H_{O_K}$, $p : \bar{C}\to \Bbb P^{d-1}_{O_{K}}$ be the projection defined by $H$, and $\bar{C}\buildrel{f}\over\to \bar{S}\to \Bbb P^{d-1}_{O_{K}}$ be the Stein factorization of $p$ (so $\bar{S}=\Spec\, p_* \CO_{\bar{C}}$), 
 $D\subset \bar{C}$ be the union of $m^{-1}(\bar{U}\smallsetminus U)$ and the exceptional divisor (viewed as a reduced scheme), $S\subset \bar{S}_K$ be the maximal open subset over which $f$ is smooth and $f|_D$ is \'etale. Set $C:= f^{-1}(S)\smallsetminus D$ and  $U' :=m(C)$; notice that $m|_C : C\iso U'$. Then $U'$, $(C,\bar{C})$, $f$ satisfy the promised properties (one has
 $f(D_f )=\bar{S}$ since  $D_f$ contains the exceptional divisor), q.e.d.
\quad\qed  \enddemo

\subhead{\rm 4.5}\endsubhead Let us return to the proof of the theorem in 4.3. We use induction by $\dim U$.
By the corollary in 4.4, we can replace $(U,\bar{U})$ by $ (C,\bar{C})$ as in loc.~cit., so we have $f:  (C,\bar{C})\to (S,\bar{S})$ with a nice section $e$ and $C$ affine over $S$. Notice that $(C,\bar{C})$ is log smooth over $(S,\bar{S})$ and the line bundle $\omega_f := \Omega_{(C,\bar{C})/(S,\bar{S})}$ equals
$ f^! (\CO_{\bar{S}} )[-1] \otimes\CO_{\bar{C}}(D_f )$. 
 
 \proclaim {Key lemma}  h-locally over $(S,\bar{S})$,  one can find a semi-stable alteration $\phi : f'\to f$ together  with a nice section $e'$ that lifts $e$ such that $(C',\bar{C}')$ is an ss-pair  and the pullback maps $\phi^* : R^1 f_* \CO_{\bar{C}} \to R^1 f'_* \CO_{\bar{C}'}$, $ f_* \omega_f \to f'_* \omega_{f'}$ are divisible by  $p$.\footnote{As elements of the groups $\Hom_{\CO_{\bar{S}}}( R^1 f_* \CO_{\bar{C}}, R^1 f'_* \CO_{\bar{C}'} )$,  $\Hom_{\CO_{\bar{S}}}( f_* \omega_f ,f'_* \omega_{f'} $).} \endproclaim
 
 For a proof, see 4.6. Assuming it for the moment, let us finish the proof of  the theorem in 4.3. By Remark (i) in 2.5, we can assume that the h-localization of $(S,\bar{S})$ in Key Lemma does not change  $\dim S$. We will show that for some
 h-covering  $(S',\bar{S}')$ of $(S,\bar{S})$ the composition $
  (C',\bar{C}')_{(S',\bar{S}')}\to 
       (C',\bar{C}')\buildrel{\phi}\over\lra
(C,\bar{C})$  is  $p$-negligible.

For any $a$ consider the exact sequence 
  $$0 \to f^* \Omega^a_{\langle S,\bar{S}\rangle} \to \Omega^a_{\langle C,\bar{C}\rangle} \to  (f^* \Omega^{a-1}_{\langle S,\bar{S}\rangle})\otimes\omega_f  \to 0. \tag 4.5.1$$ The section $e$ splits off $\Omega^a_{\langle S,\bar{S}\rangle}\hra Rf_* \Omega^a_{\langle C,\bar{C}\rangle}$ as a direct summand whose complement is $\CC one (\partial_C )$, where $\partial_C  :  \Omega^{a-1}_{\langle S,\bar{S}\rangle}\otimes f_* \omega_f  \to  \Omega^a_{\langle S,\bar{S}\rangle} \otimes R^1 f_* \CO_{\bar{C}}$,  is the boundary map for (4.5.1) (one has $R^1 f_* \omega_f =0$ since
 $f(D_f )=\bar{S}$). There is a similar splitting in case of  $f'$ provided by $e'$, and the map $\phi^* : Rf_* \Omega^a_{\langle C,\bar{C}\rangle}\to Rf'_* \Omega^a_{\langle C',\bar{C}'\rangle}$ is compatible with the direct sum decompositions. 
  Now $\phi^*$ is divisible by $p$ on the second summand: Indeed, Key Lemma asserts
 that the morphism of two-term complexes $\phi^* : \CC one (\partial_C )\to\CC one (\partial_{C'} )$ is divisible by $p$ on each term; since these are morphisms  
 of vector bundles on $O_K$-flat $\bar{S}$, our $p^{-1}\phi^*$ is uniquely defined and commutes with the differentials. Thus the map $\phi^* \otimes \id_{C_1} : \CC one (\partial_C )\otimes C_1 \to\CC one (\partial_{C'} )\otimes C_1$, where $C_1 :=\CC one (p: \Bbb Z \to \Bbb Z )$, is homotopic to zero.
  Apply $R\Gamma (\bar{S},\cdot )$ and use the induction assumption to treat the first summand $R\Gamma (\bar{S}, \Omega^a_{\langle S,\bar{S}\rangle})$;  we are done.  \quad\qed  
 
 \subhead{\rm 4.6}\endsubhead {\it Proof of Key Lemma.} Consider the relative Picard $\bar{S}$-schemes $J:= Pic^0 (\bar{C}/\bar{S})$ and $J^\flat := Pic^0 ((\bar{C}, D_f )/\bar{S})$: the first scheme parametrizes  line bundles $\CL$ on $\bar{C}$ such that the restriction of $\CL$ to the normalization of each irreducible component of any geometric fiber of $f$ has degree 0; the second one parametrizes pairs $(\CL,\gamma )$, where $\CL$ is as above and $\gamma$ is a trivialization of $\CL|_{D_f}$.
  Since $(\bar{C}, D_f )$
 is a semi-stable $\bar{S}$-family of $d$-pointed curves, $d :=\deg (D_f )$, $J$ and $J^\flat$ are semi-abelian schemes  (see \cite{R});   $J^\flat$ is an extension of $J$ by a torus $\Bbb G_m^{D_f}/\Bbb G_m$. 
  
 Over $S$ our $J^\flat$ is a  generalized Jacobian; let $i : C \to J^\flat_S$ be the  Abel-Jacobi map $i : C \to J^\flat_S$, $x\mapsto \CO_{\bar{C}} ( x-e)$. Let $C\,\,\tilde{}  \to C$ be the $i$-pullback of the multiplication by $p$ isogeny $p_{J^\flat} : J^\flat \to J^\flat$, and $ \bar{C}\,\,\tilde{} \to\bar{C}$ be the normalization of $\bar{C}$ in $C\,\,\tilde{}$. Then $f\,\,\tilde{}   :  (C\,\,\tilde{},\bar{C}\,\,\tilde{}\, )\to (S,\bar{S})$ is a family of pointed curves, which is an alteration of $f$. By the theorem in 4.4, h-locally over $(S,\bar{S})$ there is a semi-stable alteration $f'$ of $f\,\,\tilde{}$ with $(C',\bar{C}')$ semi-stable and equipped with a nice section $e'$ which lies over $e$. Let us check that the alteration $\phi : f'\to  f$ satisfies the conditions of  Key Lemma. 
  
  Set $J':=
  Pic^0 (\bar{C}'/\bar{S})$ and $J^{\prime\flat} := Pic^0 ((\bar{C}', D_f )/\bar{S})$. We have the pullback morphisms $\phi^* : J \to J' $, $J^\flat \to J^{\prime\flat}$ of our semi-abelian schemes; over $S$  we have the norm maps $\phi_{*S} : J'_S \to J_S$, $J^{\prime\flat}_S \to J^\flat_S$.  Both are  compatible with the projections $J^\flat \to J$, $J^{\prime\flat}\to J'$.  
  
     Since $\bar{S}$ is normal, for any semi-abelian $\bar{S}$-schemes $A$, $B$ one has  (see  \cite{FC} I 2.7) 
  $$\Hom (A,B)\iso \Hom(A_S ,B_S ).  \tag 4.6.1$$ Thus $\phi_{*S}$ extends to morphisms $\phi_{*} : J' \to J$, $J^{\prime\flat} \to J^\flat$.

Notice that $R^1 f_* \CO_{\bar{C}}$ is the Lie algebra of $J$, and, by Serre duality, $f_* \omega_f$ is  dual to the Lie algebra of $J^\flat$; the same is true for $f'$.
Our $ \phi^* : R^1 f_* \CO_{\bar{C}} \to R^1 f'_* \CO_{\bar{C}'}$ is the Lie algebra map for $\phi^* : J \to J'$, and $\phi^* : f_* \omega_f \to f'_* \omega_{f'}$ is the map between the duals to the Lie algebras for $\phi_* :   J^{\prime\flat} \to J^\flat$ (this is true over $S$, hence everywhere since $\bar{S}$ is $O_K$-flat). 
    
 By construction, $\phi_{*S} :  J^{\prime\flat}_S \to J^\flat_S $ factors through $p_{J^\flat}$ over $S$, i.e., it is divisible by $p$ in $\Hom (J^{\prime\flat}_S, J^\flat_S )$. By 
  (4.6.1), $\phi_*$ is divisible by $p$ in $\Hom (J^{\prime\flat}, J^\flat )$. Passing to Lie algebras, we see that  $\phi^* : f_* \omega_f \to f'_* \omega_{f'}$ is divisible by $p$. 
 Similarly, $\phi_{*S} : J'_S \to J_S$ is divisible by $p$. Notice that $J_S$, $J'_S$,
 being  Jacobians of smooth projective curves,
  are self-dual abelian schemes, and $\phi_S^* : J_S \to J'_S$ is dual to $\phi_{*S}$. Hence $\phi_S^*$ is divisible by $p$. So, by (4.6.1), $\phi^* : J \to J'$ is divisible by $p$. Passing to Lie algebras, we see that $ \phi^* : R^1 f_* \CO_{\bar{C}} \to R^1 f'_* \CO_{\bar{C}'}$ is divisible by $p$, q.e.d. \quad\qed

\head  Appendix. \endhead

Below is a proof of part (a) in the theorem from 4.4. It  follows closely
 de Jong's argument from \S\S 2--3 of \cite{dJ2} with a minor change of the lemma below; we refer the reader to sections of \cite{dJ2} for details.
 
(i) (\cite{dJ2} 2.10) One can assume that $\bar{S}$ is irreducible. By \cite{RG} 5.2.2, there is a canonical  modification of $\bar{S}$, which is projective and is an isomorphism over $S$, such that the strict transforms of $\bar{C}$ and $D_f$ are flat over $\bar{S}$. Passing to them, we can assume that {\it all fibers of $f$ have dimension 1,  of $f|_{D_f }$ have  dimension 0}.

(ii) (\cite{dJ2} 3.4--3.5) We say that a family of pointed curves is {\it good}  if irreducible components of all its geometric fibers are curves whose normalization has genus $\ge 2$. A good alteration is an alteration with good source. 

\proclaim {Lemma}  $f$ admits a good alteration h-locally over $(S,\bar{S})$.\endproclaim 

\demo{Proof of Lemma}  It suffices to find for any closed point  $s$ in $S$   its open neighborhood $S_{(s)}\subset S$ and an alteration    $(S'_{(s)},\bar{S}')$  of $(S_{(s)},\bar{S})$ such that the pullback of $f$ to  $(S'_{(s)},\bar{S}')$ admits a good alteration. To do this, we define by induction  a strictly increasing sequence of open subsets $\emptyset = V_0 \subset V_1 \subset  \ldots$
 of $\bar{S}$ and a sequence of finite extensions $F=F_0 \subset F_1 \subset \ldots$ of the field $F$ of rational functions on $\bar{C}$ such that the normalization $\bar{C}_i$ of $\bar{C}$ in $F_i$ has next properties: (a) the map $\bar{C}_i \to \bar{S}$ is smooth at $s$, (b) the map $\pi_i : \bar{C}_i \to \bar{C}$ is \'etale at $D_{f s}$,  (c) the normalizations of
irreducible components of geometric fibers of $\bar{C}_i$ over points of $V_i$ have genus $\ge 2$. There is an open neighborhood $U_i \subset S$ of $s$ over which $\bar{C}_i$ is smooth and $\pi_i$ is \'etale at $D_{f }$. The induction stops when $V_n =\bar{S}$; set $S_{(s)}=U_n$. The promised good alteration is $( \pi_n^{-1}(C ),\bar{C}_n )$ fibered over the normalization   $(S'_{(s)},\bar{S}')$   of $(S_{(s)},\bar{S})$ in $F_n$.

Let $x$ be the closed point of the closure of $s$ in $\bar{S}$. The induction produces simultaneously an auxiliary  sequence of finite subsets $T_0 \subset T_1 \subset \ldots$ of closed points of $ \bar{C}_x$; it starts with $T_0 :=$  the union of $D_{f x}$ and the set of nonregular points of $ \bar{C}_x$. The induction step:
 suppose we have  $V_{i-1}$, $F_{i-1}$, $T_{i-1}$; let us construct $V_i$, $F_i$, $T_i$
 assuming that $V_{i-1}\neq \bar{S}$.  
  Let $y$ be any closed point in $\bar{S}\smallsetminus V_{i-1}$. 
Since $\bar{S}$ is projective, there is an affine open $V$ which contains $x$ and $y$. Let $\CL$ be a very ample line bundle on $\bar{C}$. Replacing it by a sufficiently high power, we can assume that $\Gamma (\bar{C}_V ,\CL )\twoheadrightarrow \Gamma (\bar{C}_x ,\CL )\times \Gamma (\bar{C}_y ,\CL )$. 
One can find a finite unramified extension\footnote{If the residue field of $K$ is infinite, one can take $K'=K$.} $K'$ of $K$ with residue field $k'$ and two sections $\gamma_1$, $\gamma_2 \in \Gamma (\bar{C}_V ,\CL )\otimes_{O_K} O_{K'}$ which do not vanish at the generic points of irreducible components of $\bar{C}_x$, $\bar{C}_y$, such that $t=\gamma_1 /\gamma_2$ yields generically \'etale finite maps $t_x : \bar{C}_{x}\otimes k'  \to \Bbb P^1_{x}\otimes k'$, $t_y: \bar{C}_{y}\otimes k'  \to \Bbb P^1_{y}\otimes k'$  \'etale over $\{ 0,1,\infty \}$ and such that $t_x (T_{i-1}) \cap \{ 0,1,\infty \}=\emptyset$. Pick $\ell \ge 5$ prime to $p$, and let $F_i$ be an extension of $F_{i-1}$ generated by $K'$, $\mu_{\ell}$, $t^{1/\ell}$, and $(1-t)^{1/\ell}$. Let $T_i$ be the union of $T_{i-1}$ and the set of ramification points of $t_x$. The normalization $\bar{C}_i$ of $\bar{C}$ in $F_i$  satisfies (a), (b), and satisfies (c) over some open set $V_i$ which contains $V_{i-1}$ and $y$. We are done. \quad\qed  \enddemo

(iii) It remains to show that {\it every good $f$ admits   a semi-stable alteration  after a possible alteration of the base}. The genus of the generic fiber of $f$ is $\ge 2$, so  $(\bar{C}_S , D_{fS})$ is a stable  $n$-pointed curve over $S$ (where $n$ is the degree of $D_{fS}$ over $S$). The Deligne-Mumford stack of stable $n$-pointed curves is proper, so,
after replacing $(S,\bar{S})$ by an alteration, we can assume that  $ (\bar{C}_{S} , D_{fS})$ extends to a stable $n$-pointed curve $(\bar{C}' , D_{f'})$ over $\bar{S}$ (see \cite{dJ2} 3.8). We have a semi-stable family of pointed curves $f': (C',\bar{C}')\to (S,\bar{S})$, $C':= \bar{C}' \smallsetminus D_{f'}$.  By \cite{dJ2} 3.10, the goodness of $f$ implies that, after a possible alteration of $\bar{S}$, the evident morphism $\bar{C}'_{S} \to \bar{C}_S$ extends to a morphism $f'\to f$, and we are done. \qed

\enddocument